\documentclass[11pt,twoside]{article}
\usepackage{amsmath, amssymb, amsthm}
\usepackage{graphicx} 

\theoremstyle{plain}

\newtheorem{thm}{Theorem}[section]
\newtheorem{prop}[thm]{Proposition}
\newtheorem{crit}[thm]{Criterion}
\newtheorem{lemma}[thm]{Lemma}

\theoremstyle{definition}
\newtheorem{defn}[thm]{Definition}

\begin{document}

\def\Cal#1{{\cal#1}}
\def\<{\langle}\def\>{\rangle}
\def\what{\widehat}\def\wtil{\widetilde}
\def\Z{{\mathbb Z}}\def\N{{\mathbb N}} \def\C{{\mathbb C}}\def\BD{{\bf D}}
\def\Q{{\mathbb Q}}\def\R{{\mathbb R}} \def\E{{\mathbb E}}\def\BS{{\bf S}}

\def\Aut{{\rm Aut}}\def\Id{\text{Id}}
\def\lk{\text{lk}}

\def\Proof{\paragraph{Proof.}}
\def\Remark{\paragraph{Remark.}}
\def\endproof{\hfill$\square$\break\medskip}
\def\noproof{\hfill$\square$\break}
\def\noi{\noindent}

\def\notation{\paragraph{Notation.}}
\def\ackn{\paragraph{Acknowledgement.}}

\let\demph\textbf

\def\ov{\overline}

\def\al{\alpha}                 \def\be{\beta}
\def\ga{\gamma}			\def\Ga{\Gamma}\def\De{\Delta}\def\La{\Lambda}
\def\sig{\sigma}
\def\ep{\epsilon}               \def\varep{\varepsilon}
%

\title{\bf{Representations of the braid group by automorphisms of
 groups, invariants of links, and Garside groups}}
 
\author{\textsc{John Crisp and Luis Paris}}

\date{\today}

\maketitle

\begin{abstract}

From a group $H$ and a non-trivial element $h$ of $H$, we define a representation $\rho: B_n \to 
\Aut(G)$, where $B_n$ denotes the braid group on $n$ strands, and $G$ denotes the free product of 
$n$ copies of $H$. Such a representation shall be called the Artin type representation associated 
to the pair $(H,h)$. The goal of the present paper is to study different aspects of these  
representations.

Firstly, we associate to each braid $\beta$ a group $\Gamma_{(H,h)} (\beta)$ and prove that
the operator $\Gamma_{(H,h)}$ determines a group invariant of 
oriented links. We then give a topological construction of the Artin type representations 
and of the link invariant $\Gamma_{(H,h)}$, and we prove that the Artin type 
representations are faithful. The last part of the paper is dedicated to the study of some 
semidirect products $G \rtimes_\rho B_n$, where $\rho: B_n \to \Aut(G)$ is an Artin type representation.
In particular, we show that $G \rtimes_\rho B_n$ is 
a Garside group if $H$ is a Garside group and $h$ is a Garside element of $H$.
\end{abstract}



\section{Introduction}

Throughout the paper, we shall denote by $B_n$ the braid group on $n$ strands, and by $\sigma_1, 
\dots, \sigma_{n-1}$ the standard generators of $B_n$.

Let $H$ be a group, and let $h$ be a non-trivial element of $H$. Take $n$ copies $H_1, \dots, H_n$ 
of $H$, and consider the group $G=H_1 \ast \dots \ast H_n$. We denote by $\phi_i: H \to H_i$ the 
natural isomorphism and we write $h_i=\phi_i(h) \in H_i$, for all $i=1, \dots, n$.
For $k=1, \dots, n-1$, let $\tau_k: G \to G$ be the automorphism determined by
\[
\tau_k\,:\ \left\{
\begin{array}{cccl}
\phi_k(y)&\mapsto & h_k^{-1} \ \phi_{k+1}(y) \ h_k&\\
\phi_{k+1}(y)&\mapsto & h_k \ \phi_k(y)\ h_k^{-1}&\\
\phi_j(y)&\mapsto& \phi_j(y)&\text{if }\ j\neq k,k+1
\end{array}\right.
\]
for $y\in H$. One can easily show the following.

\begin{prop}\label{prop1.1}
The mapping $\sigma_k \mapsto \tau_k$, $k=1, \dots, n-1$, determines 
a representation $\rho: B_n \to \Aut(G)$.\noproof
\end{prop}

\begin{defn}
The representation of Proposition \ref{prop1.1} shall be called the {\it Artin type representation of $B_n$ 
associated to the pair $(H,h)$}.
\end{defn}

\bigskip

If $H=\Z$ and $h=1$, then $G=F_n$ is the free group of rank $n$ and $\rho$ is the classical 
representation introduced by Artin in \cite{Art1}, \cite{Art2}. Another example which appears in the 
literature is the case where $H=\Z$ and $h$ is a non-zero integer. This case has been introduced 
by Wada \cite{Wad} in his construction of group invariants of links. Sections 2 and 3 of the 
present paper are inspired by Wada's work \cite{Wad}. 

\bigskip
Our purpose in this paper is to study different aspects of the Artin type representations.

\begin{defn}
Let $\rho: B_n \to \Aut(G)$ be the Artin type representation associated to a 
pair $(H,h)$. Let $\beta \in B_n$. Then we denote by $\Gamma(\beta) = \Gamma_{(H,h)} (\beta)$ the 
quotient of $G$ by the relations
\[
g = \rho (\beta) g, \quad g \in G.
\]
\end{defn}
\medskip

For a braid $\beta$, we denote by $\what \beta$ the oriented link (or more precisely the
equivalence class of oriented links) represented by the closed braid of $\beta$ as defined in \cite{Bir}. 
 Given two braids $\beta_1$ and $\beta_2$ (not necessarily with 
the same number of strands), we prove in Section \ref{Sect2} that $\Gamma (\beta_1) \simeq \Gamma 
(\beta_2)$ if $\what \beta_1 = \what \beta_2$. This allows us to define a group invariant of oriented links, 
$\Gamma_{(H,h)}$, by setting $\Gamma_{(H,h)} (L)$ to be the group $\Gamma_{(H,h)} (\beta)$ for any 
braid $\beta$ such that $L=\what\beta$.
Note that, in the case $H=\Z$ and $h=1$, the invariant $\Gamma_{(\Z,1)}$ computes the link group,
namely  $\Gamma_{(\Z,1)}(L)\cong\pi_1(\BS^3\setminus L)$ for any link $L$ in $\BS^3$.

The goal of Section \ref{Sect3} is to give topological constructions of the Artin type representations and 
of the groups $\Gamma_{(H,h)} (\beta)$, for $\beta \in B_n$. If $H=\Z$ and $h$ is a non-zero integer, 
then our constructions coincide with Wada's constructions (see \cite{Wad}, Section 3). In fact, 
our constructions are straightforward extensions of Wada's constructions to all Artin type 
representations.

In Section \ref{Sect4}, we prove that Artin type representations are faithful (Proposition \ref{4.1}).
If $h$ has  infinite order, then the Artin type representation $\rho: B_n \to \Aut(G)$
contains the classical Artin representation and, therefore,
is faithful by \cite{Art1}, \cite{Art2}. So, Proposition \ref{4.1} is mostly of 
interest in the case where $h$ has finite order. In fact the proof may be easily reduced
to the case $H= \Z /k\Z$ and $h=1$, however we will not need to use any such reduction, as
our method applies just as easily in all cases. We note also that the case where $H$ is 
cyclic of order $2$ follows (by somewhat different methods) from Section 2.3 of \cite{CP}.
The proof of  Proposition \ref{4.1} is inspired by the proof of Theorem A of \cite{Shp},
and it is based on Dehornoy's work on orders on braids \cite{Deh1}, \cite{Deh2}.

The remaining sections (Sections \ref{Sect5} and \ref{Sect6}) are
dedicated to the study of semidirect products $G \rtimes_\rho B_n$, where $\rho: B_n 
\to \Aut(G)$ is the Artin type representation associated to a pair $(H,h)$.

If $H=\Z$ and $h=1$, then $G \rtimes_\rho B_n$ is the Artin group $A(B_n)$ associated to the Coxeter 
graph $B_n$ (not to be confused with the braid group $B_n$, which is itself an
Artin group, of type $A_{n-1}$). This result is implicit in \cite{Lam}, \cite{Cri},
and explicit in \cite{CP}.
The group $A(B_n)$ is well-understood. In particular, solutions to the word and conjugacy problems in this 
group are known (see \cite{Del}, \cite{BS}), it is torsion free (see \cite{Bri}, \cite{Del}),
its center is an infinite cyclic group (see \cite{Del}, \cite{BS}),
it is biautomatic (see \cite{Cha1}, \cite{Cha2}), and it has an 
explicit finite dimensional classifying space (see \cite{Del}, \cite{Bes}).

A natural next step is to understand the groups $G \rtimes_\rho B_n$ in 
the case where $\rho$ is a Wada representation (of type 4), namely, when $H=\Z$ and $h\in \Z 
\setminus \{0\}$. One can readily establish that, for these representations,
the group $G \rtimes_\rho B_n$ fails to be an Artin group unless $h=\pm 1$.
It turns out, however, that these groups do have quite a lot in common with Artin groups: like
the Artin groups, they  belong to a family of groups known as \emph{Garside groups}.

Briefly, a \emph{Garside group} is a group $G$ which admits a left invariant lattice order and contains
a so-called \emph{Garside element}, a positive element $\Delta$ whose positive divisors generate $G$
and such that conjugation by $\Delta$ leaves the lattice structure invariant (there are also conditions
placed on the positive cone of $G$, that it be a finitely generated atomic monoid -- see Section \ref{Sect5} for
more details).
The notion of a Garside group was introduced by Dehornoy and the second author \cite{DP}
in a slightly restricted sense, and, later, by Dehornoy \cite{Deh5} in the larger sense which is
now generally used. The theory of Garside groups is largely inspired by
the papers of Garside \cite{Gar}, which treated the case of braid groups, and Brieskorn and Saito \cite{BS}
which generalised Garside's work to Artin groups. 
The Artin groups of spherical (or finite) type which include, 
notably, the braid groups as well as the groups $A(B_n)$ mentioned above,
are motivating examples. 
Other interesting examples of Garside groups
include all torus link groups (see \cite{Pic3}) and some generalized braid groups associated
to finite complex reflection groups. 

Garside groups have many attractive properties.
Solutions to the word and conjugacy problems in these groups are 
known (see \cite{Deh5}, \cite{Pic1}, \cite{FG}), they are torsion free (see \cite{Deh4}),
they admit canonical 
decompositions as iterated crossed products of ``irreducible'' components, and the center of each 
component is an infinite cyclic group (see \cite{Pic2}), they are biautomatic (see \cite{Deh5}), and they 
admit finite dimensional classifying spaces (see \cite{DL}, \cite{CMW}). Another important property of the 
Garside groups is that there exist criteria in terms of presentations to detect them (see 
\cite{DP}, \cite{Deh5}).

In Section \ref{Sect6}, we prove that, if $H$ is a Garside group, $h$ a Garside element of $H$,
and $\rho$ the Artin type representation associated to $(H,h)$, then $G\rtimes_\rho B_n$
is also a Garside group (Theorem \ref{6.1}). This result applies in particular to 
the case $H=\Z$ and $h \in \Z \setminus \{ 0 \}$, but also applies, for example, to the case 
where $H$ is another braid group, say $H=B_l$, and $h=\Delta^k$ is a non-trivial power of the 
fundamental element of $B_l$.

The proof of Theorem \ref{6.1} is based on a criterion which is developed in Section \ref{Sect5}
for proving that a given group is a Garside group. This criterion rests largely on the
``coherence" condition of \cite{DP} which has its roots in the original arguments of Garside \cite{Gar}.
It is essentially a variation of other criteria appearing in the literature
(see, for example, \cite{Deh5} Prop. 6.14). Our criterion differs from that of Dehornoy \cite{Deh5}
just mentioned in that it is not algorithmic. In particular, we do not give any method for
finding a Garside element. However, it is relatively easy to apply once one has an appropriate
presentation and an expression for a Garside element to hand.  
 
\medskip

Finally, we add an appendix to our paper in order to answer a question posed by Shpilrain in
his study of Wada's representations \cite{Shp}, and which is otherwise a little tangential to 
the main subject of this paper.

\begin{defn}
Let $G$ be a group. Two 
representations $\rho, \rho': B_n \to \Aut( G)$ are called {\it equivalent} if there exist 
automorphisms $\phi: G \to G$ and $\mu: B_n \to B_n$ such that
$$
\rho' (\mu( \beta)) = \phi^{-1} \circ \rho( \beta) \circ \phi
$$
for all $\beta \in B_n$.
\end{defn}

\bigskip\noindent
{\bf Remark.} 
If two representations $\rho, \rho': B_n \to \Aut(G)$ are equivalent, then the groups $G 
\rtimes_\rho B_n$ and $G \rtimes_{\rho'} B_n$ are isomorphic.

\bigskip
Shpilrain's question (see \cite{Shp}) was simply to give a classification of Wada's 
representations up to equivalence. This classification is given in Proposition A.1.


\section{Link invariants}\label{Sect2}

Let $H$ be a group, $h$ a non-trivial element of $H$, and $\rho: B_n \to \Aut(G)$ be the Artin 
type representation associated to $(H,h)$. Recall that the group $G$ is defined as
$G= H_1 \ast \dots \ast H_n$, where group isomorphisms $\phi_i:H_i\to H$ are given for $i=1,2,\dots ,n$.
The goal of this section is to prove the following.

\begin{prop}\label{2.1} 
Let $n,m \in \N$, and let 
$\beta_1 \in B_n$ and $\beta_2 \in B_m$. If $\what \beta_1 = \what \beta_2$, then 
$\Gamma_{(H,h)}(\beta_1) \simeq \Gamma_{(H,h)}(\beta_2)$.
\end{prop}

\begin{defn}
Let $L$ be an oriented link. We set  $\Gamma_{(H,h)}(L) := \Gamma_{(H,h)}(\beta)$,
where $\beta$ is any braid (on any number of strings)
such that $L={\what\beta}$. By Proposition \ref{2.1}, $\Gamma_{(H,h)}$ is 
a well-defined group invariant of oriented links.
\end{defn}

\paragraph{Proof of Proposition \ref{2.1}.}
Let $n \in \N$ and let $\beta \in B_n$. We write $\Gamma$ for $\Gamma_{(H,h)}$.
By Markov's theorem (see [Birm, Thm. 2.3]), it suffices to show:

\smallskip
(1) $\Gamma (\alpha^{-1} \beta \alpha) \simeq \Gamma (\beta)$ for all $\alpha \in B_n$;

\smallskip
(2) $\Gamma(\beta \sigma_n) \simeq \Gamma(\beta)$;

\smallskip
(3) $\Gamma (\beta \sigma_n^{-1}) \simeq \Gamma (\beta)$;

\smallskip\noindent
where $\beta \sigma_n$ and $\beta \sigma_n^{-1}$ are viewed as braids on $n+1$ strands.
\smallskip

For a given $n\in\N$, we write $G_{(n)}=H_1\ast\dots\ast H_n$. Note that, if $\beta\in B_n$ and
$n\leq m$, then  the action of $\beta$ via $\rho$
on $G_{(m)}$ agrees with the action via $\rho$ on $G_{(n)}<G_{(m)}$, and is trivial on the free
factors $H_{n+1},..,H_m$. We suppress $\rho$ from our notation, writing simply $\beta(g)$ to mean
$\rho(\beta)g$, for any $\beta\in B_n$ and $g\in G_{(m)}$.

\paragraph{\it Proof of (1):}
For $\beta\in B_n$, the
group $\Gamma(\beta)$ is defined as the quotient of $G_{(n)}$ by the relations 
$g=\beta(g)$ for all $g\in G_{(n)}$. Since, for $\alpha\in B_n$,
 the relation $g=\alpha^{-1}\beta\alpha(g)$ is equivalent to the relation
$\alpha(g) =\beta(\alpha(g))$, and $\alpha$ is an automorphism of $G_{(n)}$, it is clear 
that $\Gamma(\alpha^{-1}\beta\alpha)$ is defined by the same set of relations as $\Gamma(\beta)$.
So (1) holds.
 
\paragraph{\it Proof of (2):}
The group $\Gamma(\beta\sigma_n)$ may be defined as the quotient of $G_{(n+1)}$ by
the family of relations $R(i,x):\phi_i(x)= \beta\sigma_n(\phi_i(x))$ for all $i=1,2,..,n+1$ and
all $x\in H$. Note that  
$\sigma_n(\phi_{n+1}(x))=h_n\phi_{n}(x)h_n^{-1}$. Therefore the relation $R(n+1,x)$
is equivalent to the relation $R'(n+1,x): \phi_{n+1}(x)=\beta(h_n\phi_{n}(x)h_n^{-1})$, where
the right hand side is actually an element of $G_{(n)}$. In particular $\Gamma(\beta\sigma_n)$
is generated by the image of $G_{(n)}$. Also, 
$$
\beta\sigma_n(\phi_n(x))=\beta(h_n^{-1}\phi_{n+1}(x)h_n)=\beta(h_n^{-1})\phi_{n+1}(x)\beta(h_n)\,.
$$
So, in view of $R'(n+1,x)$, the relation $R(n,x)$ is now equivalent to the relation
$R'(n,x): \phi_n(x) = \beta(\phi_n(x))$. Finally,
since $\sigma_n(\phi_i(x))=\phi_i(x)$ for all $i<n$, the remaining relations $R(i,x)$ are
equivalent to $R'(i,x):\phi_i(x)=\beta(\phi_i(x))$ for
all $i=1,2,..,n-1$, and all $x\in H$.
It now follows that $\Gamma(\beta\sigma_n) \simeq \Gamma(\beta)$.

\paragraph{\it Proof of (3):} 
Observe that $\Gamma(\beta^{-1})\simeq\Gamma(\beta)$, since the relation $g=\beta(g)$
is equivalent to $\beta^{-1}(g)=g$, for all $g\in G_{(n)}$. Then 
\[
\begin{array}{rll}
\Gamma(\beta\sigma_n^{-1})&\simeq \Gamma(\sigma_n\beta^{-1})&\\
&\simeq\Gamma(\beta^{-1}\sigma_n)\hskip3mm&\text{by part (1) of the proof,}\\
&\simeq\Gamma(\beta^{-1})&\text{by part (2) of the proof,}\\
&\simeq\Gamma(\beta)\,.&
\end{array}
\]\endproof



\section{Topological construction of the link invariants}\label{Sect3}

\bigskip
Let $X$ be a CW-complex, let $P_0 \in X$ be a basepoint, and let $\alpha: [0,1] \to X$ be a loop 
based on $P_0$. We assume that $\alpha$ is not homotopic to the constant path.
In this section we give a topological realization of the Artin type representation of $B_n$
associated to the pair $(H,h)= (\pi_1 (X, P_0), [\alpha])$, and we deduce a topological construction
of the link invariant $\Gamma_{(H,h)}$ of the previous section.

Let $\BD= \BD({n+1 \over 2}, {n+1 \over 2})$ denote the disk in $\C$ centered at ${n+1 \over 2}$ 
of radius ${n+1 \over 2}$. Now, we construct a space $Y$ obtained from $\BD$ by making $n$ holes in 
$\BD$ and gluing a copy of $X$ into each hole by identifying the circular boundary of the hole to the 
loop $\alpha$ in $X$.

Choose some small $\varepsilon >0$ (we require only that $\varepsilon < {1 \over 8}$). Let
$$
Y'=\BD \setminus \left(\bigcup_{k=1}^n \BD^\circ (k,\varepsilon) \right),
$$
where $\mathop{\BD}^\circ (k,\varepsilon)$ denotes the open disk centered at $k$ of radius 
$\varepsilon$. Take $n$ copies $X_1, \dots, X_n$ of $X$, denote by $f_k: X \to X_k$ the natural 
homeomorphism, and write $\alpha_k = f_k \circ \alpha$ for all $k=1, \dots, n$. Then
$$
Y= \left( Y' \sqcup \left( \bigsqcup_{k=1}^n X_k \right)\right) / \sim,
$$
where $\sim$ is the identification defined by
$$
\alpha_k(t) \sim k+ \varepsilon e^{2i\pi t},\quad k=1, \dots, n,\ t\in [0,1].
$$

Finally, choose a basepoint $Q_0\in\partial\BD$ for $Y$.
The following result is a direct consequence of the above construction.

\begin{lemma}\label{3.1} 
Let $H=\pi_1(X,P_0)$, and let $H_1, \dots, H_n$ be $n$ copies of $H$. Then 
$\pi_1(Y,Q_0)\simeq H_1 \ast \dots \ast H_n$.\noproof
\end{lemma}

We now show that the braid group $B_n$ acts on $Y$ up to isotopy relative to the boundary of 
$\BD$ in such a way that the induced action on $\pi_1(Y)$ is the Artin type representation
associated to  $(H,h)$, where $h$ is the element of $H=\pi_1(X,P_0)$ represented by $\alpha$.

\begin{figure}[ht]
\begin{center}
\includegraphics[width=10cm]{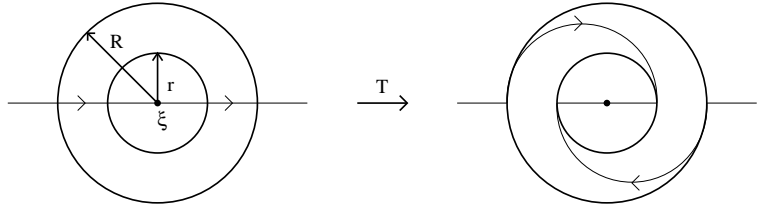}
\end{center}
\caption{A half Dehn twist.}\label{figDehntwist}
\end{figure}

Let $\xi \in \C$ and $0<r<R$. Define the {\it half Dehn twist} $T=T(\xi,r,R)$ by
\[
T(\xi+\rho e^{i \theta}) = \begin{cases}
\xi+\rho e^{i( \theta -\pi)} &\text{if }\, 0\le \rho \le r\\
\xi+\rho e^{i( \theta -t\pi)} &\text{if }\, r \le \rho \le R\text{ and } t=\frac{R-\rho}{R-r}\\
\xi +\rho e^{i\theta} &\text{if }\, \rho \ge R
\end{cases}
\]
(see Figure \ref{figDehntwist}).

\begin{figure}[ht]
\begin{center}
\includegraphics[width=12cm]{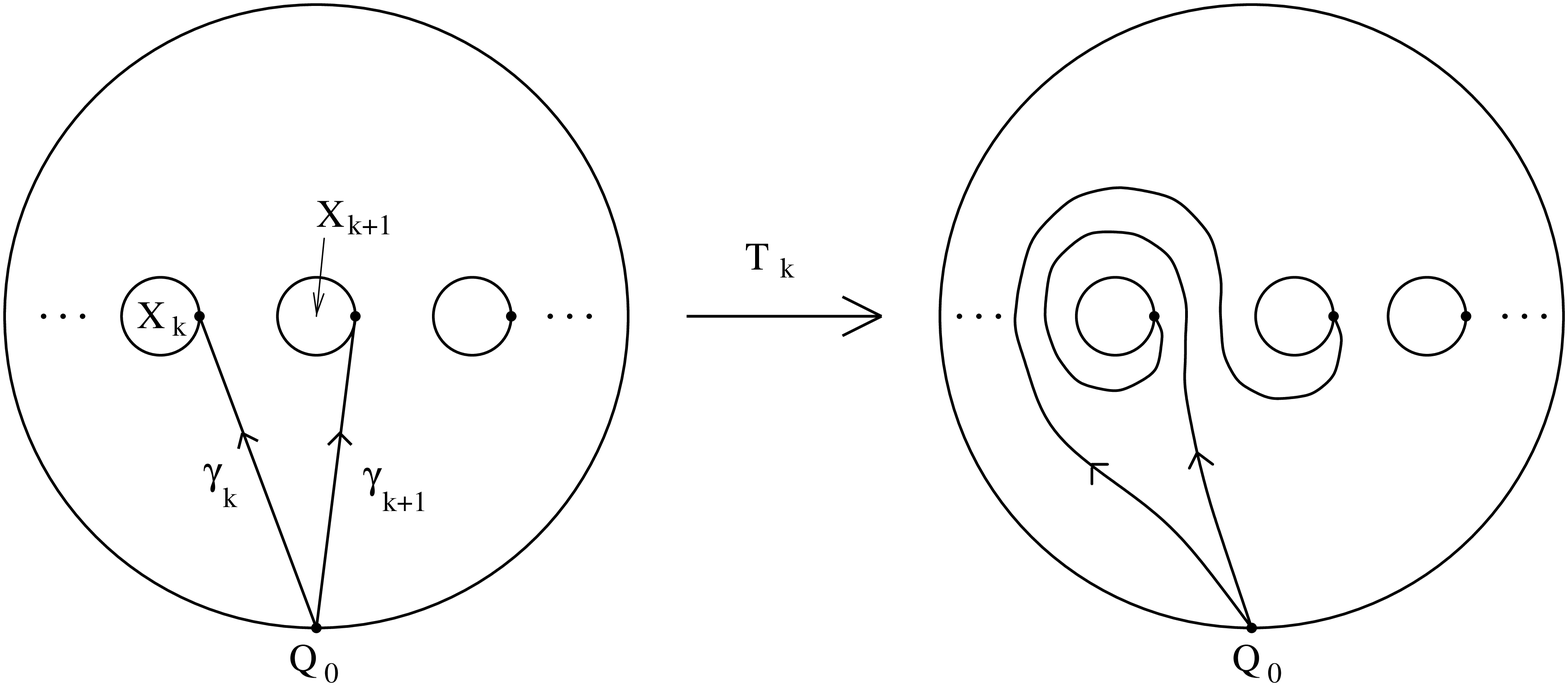}
\end{center}
\caption{The homeomorphism $T_k:Y'\to Y'$.}\label{figTk}
\end{figure}

\bigskip
Let $T_k^\BD: \BD \to \BD$ be the homeomorphism defined by
$$
T_k^\BD= T(k,\varepsilon, 2\varepsilon)^{-3} \circ T(k+1, \varepsilon, 2\varepsilon)^{-1} \circ T(k+{1 \over 
2}, {1 \over 2} +\varepsilon, {1 \over 2} +2\varepsilon).
$$
Note that $T_k^\BD$ leaves invariant the set $\cup_{j=1}^n \BD (j,\varepsilon)$, and therefore restricts
to a homeomorphism $T_k: Y' \to Y'$. See Figure \ref{figTk}.

One can verify (with a little effort) that $T_k T_{k+1} T_k$ is isotopic to $T_{k+1} T_k T_{k+1}$
relative to $\partial Y'$ for $k=1, \dots, n-2$, and that $T_k T_l$ is 
isotopic to $T_l T_k$ relative to $\partial Y'$ for $|k-l| \ge 2$. Moreover, $T_k$ fixes $\partial\BD$
and transforms the rest of $\partial Y'$ as follows:
\[
T_k(j+\varepsilon e^{i\theta})=\begin{cases}
j+\varepsilon e^{i\theta} &\text{if }\, j\neq k,k+1\\
k+1+\varepsilon e^{i\theta} &\text{if }\,j=k\\
k+\varepsilon e^{i\theta} &\text{if }\, j=k+1\,.
\end{cases}
\]
Therefore, $T_k$ extends to a homeomorphism $T_k: Y \to Y$ by setting, for all $x\in X$,
\[
T_k(f_j(x))=\begin{cases}
f_j(x) &\text{if }\, j\neq k,k+1\\
f_{k+1}(x) &\text{if }\, j=k\\
f_k(x) &\text{if }\, j=k+1\,.
\end{cases}
\]
The homeomorphism $T_k$ is the identity on $\partial \BD$, $T_k T_{k+1} T_k$ 
is isotopic to $T_{k+1} T_k T_{k+1}$ relatively to $\partial \BD$ for $k=1, \dots, n-2$, and $T_k 
T_l$ is isotopic to $T_l T_k$ relatively to $\partial\BD$ for $|k-l| \ge 2$.

 By the above observations, $T_k$ determines an 
automorphism $\tau_k: \pi_1(Y,Q_0) \to \pi_1(Y,Q_0)$. Moreover,
\[
\begin{array}{rll}
\tau_k\tau_{k+1}\tau_k\! &= \tau_{k+1} \tau_k \tau_{k+1}\hskip5mm &{\rm for}\ k=1, \dots, n-2\\
\tau_k \tau_l &= \tau_l \tau_k\! &{\rm for}\ |k-l|\ge 2\\
\end{array}
\]

\noindent
So, the mapping $\sigma_k \to \tau_k$ determines a representation $\rho: B_n \to \Aut( 
\pi_1(Y,Q_0))$.

\begin{figure}[ht]
\begin{center}
\includegraphics[width=5cm]{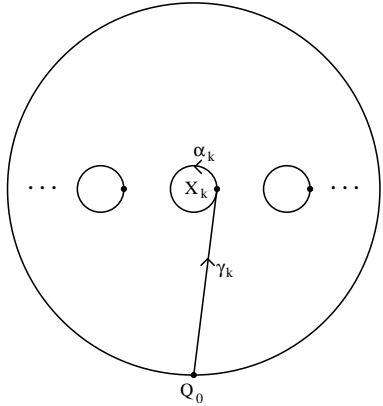}
\end{center}
\caption{The path $\gamma_k$.}\label{figGammak}
\end{figure}

Set $Q_0= {n+1\over 2} +i{n+1 \over 2}$. Let $\gamma_k:[0,1] \to Y$ denote the path which joins 
$Q_0$ to $f_k(P_0)$ represented in Figure \ref{figGammak}. We identify $\pi_1(Y,Q_0)$ with
$G=H_1\ast\dots\ast H_n$ in such a way that the $k$-th embedding $\phi_k: 
H=\pi_1(X,P_0) \to H_k \subset G$ is defined by
$$
\phi_k([\beta]) = [\gamma_k \beta \gamma_k^{-1}].
$$
With this assumption, one can easily show the following.

\begin{prop}\label{3.2}
The representation $\rho: B_n \to \Aut( \pi_1(Y,Q_0))$ described 
above coincides with the Artin type representation of $B_n$ associated to $(H,h)$, where 
$H=\pi_1(X,P_0)$ and $h$ is the element of $H$ represented by $\alpha$.\noproof
\end{prop}

Consider an oriented $m$-component link $L=K_1 \cup \dots \cup K_m$ in $\BS^3$. 
The knot $K_i$ is an embedding $K_i: \BS^1 \to \BS^3$, and $K_i(\BS^1) \cap K_j(\BS^1) = 
\emptyset$ for $i \neq j$. Define a {\it tubular neighborhood} of $K_i$ to be an embedding $T_i: 
\BD^2 \times \BS^1 \to \BS^3$ such that $T_i(0,\xi)= K_i(\xi)$ for all $\xi \in \BS^1$. Here, 
$\BD^2$ denotes the disk centered at $0$ of radius 1 in $\C$. A {\it framing} of $L$ is a 
collection $\{ T_i: \BD^2 \times \BS^1 \to \BS^3\}_{i=1}^m$ of embeddings such that $T_i$ is a 
tubular neighborhood of $K_i$, for $i=1, \dots ,m$, and $T_i( \BD^2 \times \BS^1) \cap T_j( \BD^2 
\times \BS^1) = \emptyset$ for $i \neq j$. The {\it longitude} of the component
$K_i$ is the (oriented) embedding $\lambda_i: \BS^1\to\BS^3$ such that $\lambda_i(\xi)=T_i(1,\xi)$
 for all $\xi \in \BS^1$. The tubular neighborhood of the framing of each component 
$K_i$ is determined up to isotopy
by the homology class of its longitude $\lambda_i$ in the knot complement $\BS^3\setminus K_i$.

\begin{figure}[ht]
\begin{center}
\includegraphics[width=10cm]{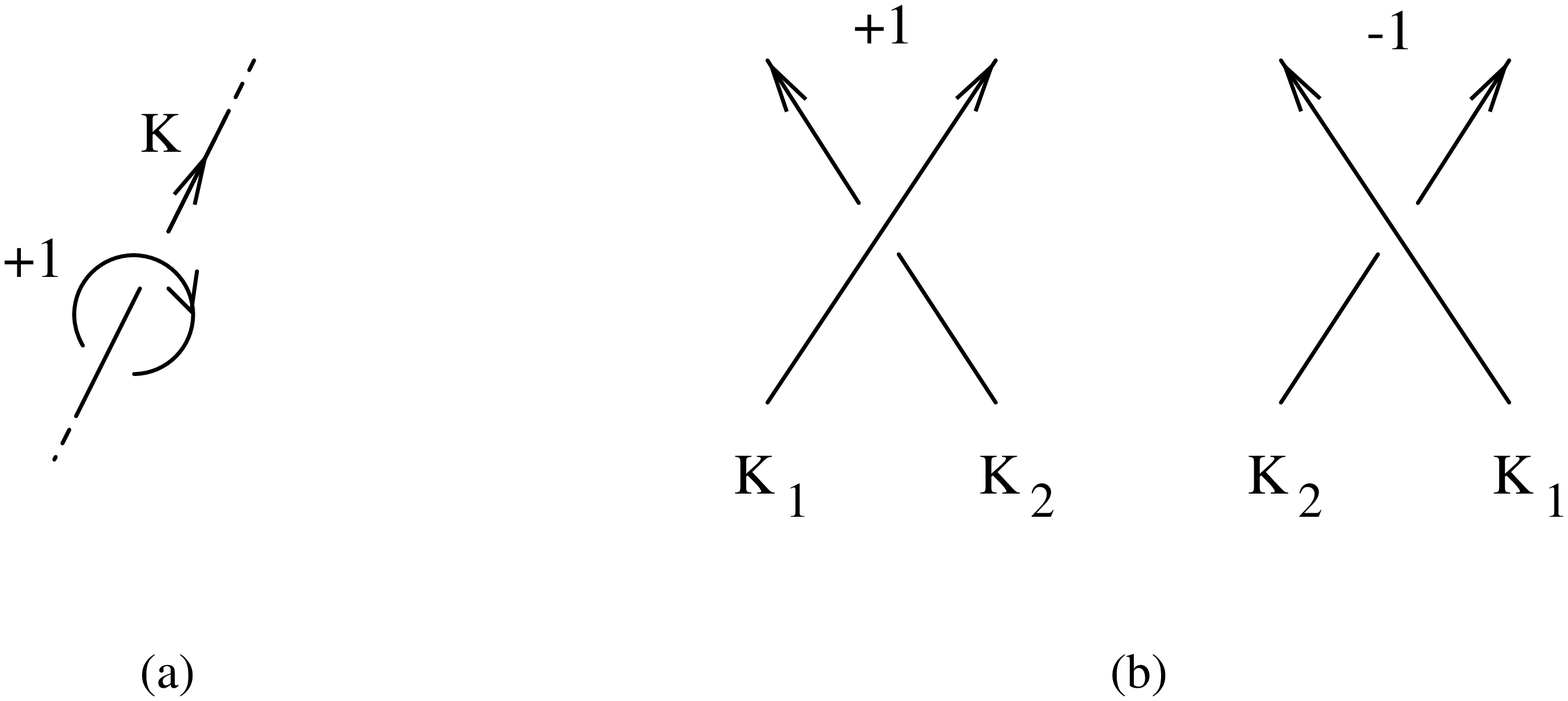}
\end{center}
\caption{Sign conventions.}\label{figSigns}
\end{figure}

Given an oriented knot $K$, we identify $H_1(K):=H_1(\BS^3\setminus K)$ with $\Z$
in such a way that $1\in\Z$ is represented by the 1-cycle depicted in Figure \ref{figSigns}(a).
Let $K_1$, $K_2$ denote disjoint oriented knots in $\BS^3$. One defines the \emph{linking number}
$\lk(K_1,K_2)\in \Z$ to be the class $[K_1]\in H_1(K_2)=\Z$. The linking number $\lk(K_1,K_2)$ may be
measured from any regular projection of the link $K_1\cup K_2$ by counting with sign the crossings
where $K_1$ passes over $K_2$, as indicated in Figure \ref{figSigns}(b). (Equally one may choose to count 
undercrossings with the appropriate sign, and one quickly sees that $\lk(K_1,K_2)=\lk(K_2,K_1)$).

\paragraph{Notation (Preferred framing).} 
Let $L=K_1 \cup \dots \cup K_m$ be an $m$-component oriented link in $\BS^3$. 
Up to isotopy, there is a unique framing in which the longitude $\lambda_i$ for each component $K_i$
satisfies the following condition:
\[
\sum_{j=1}^m \lk(\lambda_i, K_j)=0\,.
\]
Note that, for $j\neq i$, $\lk(\lambda_i,K_j)=\lk(K_i,K_j)$ and is determined by the oriented link $L$.
We shall refer to the above framing as the \emph{preferred framing} of $L$.

We now wish to associate to an oriented link $L$ the space $\Omega(L,X)$ obtained by performing a
`generalised' surgery on the link $L$ according to the preferred framing just described.
More precisely, let $L=K_1 \cup \dots \cup K_m$ and let $\{ T_i: \BD^2 \times \BS^1 \to \BS^3\}_{i=1}^m$
be the preferred framing. 
Let $\mathop{T_i}^\circ$ denote the interior of $T_i(\BD^2 \times \BS^1)$ for $i=1,\dots ,m$, 
and let 
$$
\Omega'(L)= \BS^3 \setminus \left( \bigcup_{i=1}^m T_i^\circ \right).
$$
Take $m$ copies $X_1, \dots, X_m$ of $X$, denote by $f_i: X \to X_i$ the natural homeomorphism, 
and write $\alpha_i= f_i \circ \alpha$. Then
$$
\Omega(L,X) =\left( \Omega'(L) \sqcup\left( \bigsqcup_{i=1}^m (X_i \times \BS^1) \right)\right) / \sim,
$$
where $\sim$ is the identification defined by
$$
(\alpha_j(t), \eta) \sim T_j(e^{2i\pi t}, \eta), \quad j=1, \dots, m,\ t\in [0,1],\ \eta \in \BS^1.
$$

The following proposition yields a second proof of the fact that $\Gamma_{(H,h)}$ is a link invariant
for any finitely generated group $H$ and nontrivial element $h\in H$. 

\begin{prop}\label{3.3}
Let $\beta$ be a braid, and let $\what\beta$ denote the closed braid of $\beta$. Let $X$ be a CW-complex
with basepoint $P_0$ and let $\alpha$ be a nontrivial loop in $X$.
Then $\pi_1( \Omega( \what \beta,X))$ is isomorphic to $\Gamma_{(H,h)}(\beta)$,
where $H=\pi_1(X,P_0)$ and $h$ is the element of $H$ represented by $\alpha$.
\end{prop}

\Proof 
We first remind the reader of the standard construction of the closed braid $\what\beta$ from
a braid $\beta$ (see [Birm]). Firstly, decompose $\BS^3$ as follows: let $T_1,T_2$ be two
copies of the solid torus $\BD\times\BS^1$ and write
\[
\BS^3=T_1\bigcup_{\kappa:\partial T_1\to \partial T_2} T_2\,,
\]
where the identifying map $\kappa$ is a homeomorphism carrying $\partial \BD$ to $\BS^1$ and $\BS^1$ to 
$\partial\BD$. 
The closed braid $\what\beta$ is the oriented link which is induced by composing the braid
$\beta: \{1,\dots ,n\}\times [0,1]\to \BD\times [0,1]$ with 
the composition of maps
$$
\BD\times[0,1]\buildrel f\over\longrightarrow \BD\times\BS^1=T_1\buildrel g\over\longrightarrow\BS^3\,,
$$  
where $f(p,t)=(p,e^{i2\pi t})$, for $p\in\BD$ and $t\in [0,1]$,
and $g$ denotes the inclusion of $T_1$ in $\BS^3$. The orientation on $\what\beta$ is naturally
induced from a choice of orientation of the interval $[0,1]$.

\begin{figure}[ht]
\begin{center}
\includegraphics[width=10cm]{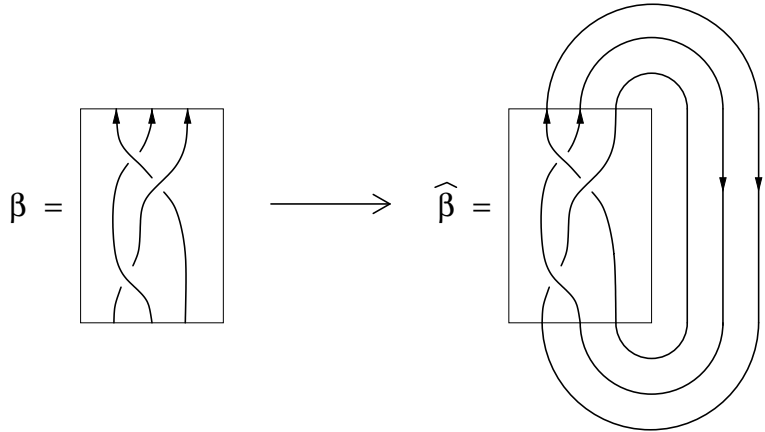}
\end{center}
\caption{Braid closure.}\label{figClosure}
\end{figure}

\begin{figure}[ht]
\begin{center}
\includegraphics[width=6cm]{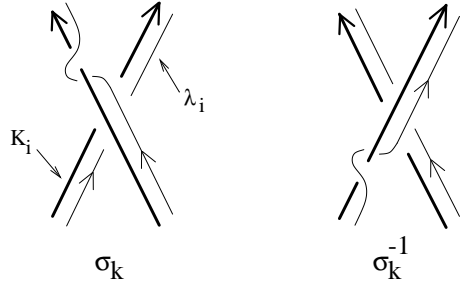}
\end{center}
\caption{Choosing a framing for $\what\beta=K_1\cup\dots\cup K_m$.}\label{figFraming}
\end{figure}

 Given a standard projection of a braid $\beta$ 
we may describe a projection of the closed braid $\what\beta$ with the same number of crossings,
as indicated in Figure \ref{figClosure}. We now produce a framing $\Lambda$ of $\what\beta$ by choosing
a longitude $\lambda_i$ for each component $K_i$ of $\what\beta$ whose projections are as indicated in
Figure \ref{figFraming} in the vicinity of a crossing, and otherwise parallel to the link projection.
It is easily enough verified that this framing is exactly the preferred framing of $\what\beta$.

Write $\be=\sig_{i_1}^{\varepsilon_1}\sig_{i_2}^{\varepsilon_2}\dots\sig_{i_r}^{\varepsilon_r}$ 
and define $T^\BD_\be:\BD\to \BD$ as the composition
of the homeomorphisms $(T_{i_j}^\BD)^{\varepsilon_j}$ for $j=1, \dots, r$.
Similarly, define $T_\beta = T_{i_1}^{\varepsilon_1} T_{i_2}^{\varepsilon_2} \dots T_{i_r}^{\varepsilon_r}:
Y \to Y$.
For $j=1,..,n$, denote by $b_j$ the point $j+\varep$ on $\partial \BD(j,\varep)$. 
This is the point on $\partial Y'$ to which the basepoint of $X_j$ is attached when forming $Y$. Since
$T_\be^\BD$ is isotopic to $\Id_\BD$ relative to $\partial\BD$, there is a homeomorphism 
$U:\BD\times [0,1]\to \BD\times [0,1]$ such that
$U(x,0)=(x,0)$, $U(x,1)=(T_\be^\BD(x),1)$, for all $x\in\BD$,
and $U$ fixes $\partial \BD \times [0,1]$ pointwise.  Moreover,
by construction, $U$ carries $(\bigsqcup\limits_{j=1}^n \BD(j,\varep))\times [0,1]$
to a tubular neighbourhood of (a representative of) the braid $\be$, and $g \circ f \circ U$ carries the arcs
$\{ b_j\times [0,1] : j=1,..,n \}$ to a framing of $\what\be$ equivalent to that
described in Figure \ref{figFraming}, namely the preferred framing.
Consequently the space $\Omega(\what\be, X)$ is homeomorphic to $T'_1\cup T_2$ where
\[
T_1'=Y\times [0,1]/((y,0)\sim (T_\be(y),1))\,.
\]
We therefore have $\pi_1(T_1')\cong G\star \< t\> /(tgt^{-1}\sim \rho(\be)g\ \forall\, g\in G)$, an HNN-extension.
Attaching $T_2$ to $T_1'$ has the effect of  simply killing the stable letter $t$. Consequently
\[
\pi_1(\Omega(\what\be, X))\cong G/(g\sim \rho(\be)g\ \forall\, g\in G)= \Gamma_{(H,h)}(\beta)\,.
\] 
\endproof


\section{Faithfulness}\label{Sect4}

\bigskip
Consider a group $H$ and a non-trivial element $h \in H \setminus \{1\}$, and write $G=H_1 \ast 
\dots \ast H_n$, where $H_i$ is a copy of $H$. The aim of this section is to prove the following.

\begin{prop}\label{4.1}
Let $\rho: B_n \to \Aut(G)$ be the Artin type representation of $B_n$ 
associated to $(H,h)$. Then $\rho$ is faithful.
\end{prop}

As pointed out in the introduction,
the proof of Proposition \ref{4.1} is strongly inspired by the proof of Theorem A of [Shp], and its 
main ingredient is the following result due to Dehornoy [Deh1], [Deh2].

\begin{prop}[Dehornoy]\label{4.2}
Let $B_{n-1}$ be the subgroup of $B_n$ generated by 
$\sigma_2, \dots, \sigma_{n-1}$. Let $\beta \in B_n$. Then either

\smallskip
(1) $\beta \in B_{n-1}$; or

\smallskip
(2) $\beta$ can be written 
$$
\beta= \alpha_0 \sigma_1 \alpha_1 \sigma_1\alpha_2 \dots \sigma_1 \alpha_l,
$$
where $l\ge 1$ and $\alpha_0, \dots, \alpha_l \in B_{n-1}$; or

\smallskip
(3) $\beta$ can be written 
$$
\beta= \alpha_0 \sigma_1^{-1} \alpha_1 \sigma_1^{-1}\alpha_2 \dots \sigma_1^{-1} \alpha_l,
$$
where $l\ge 1$ and $\alpha_0, \dots, \alpha_l \in B_{n-1}$.\noproof
\end{prop}

The following lemma is a preliminary result to the proof of Proposition \ref{4.1}.

\begin{lemma}\label{4.3}
Let $G'=H_2 \ast \dots \ast H_n$. Let $u\in G$ such that the normal form of 
$u$ with respect to the decomposition $G=H_1 \ast G'$ starts with $h_1^{-1}$ and ends with $h_1$.

\smallskip
(1) The normal form of $\rho(\sigma_1)(u)$ with respect to the decomposition $G=H_1 \ast G'$ 
also starts with $h_1^{-1}$ and ends with $h_1$.

\smallskip
(2) Let $k \in \{2, \dots, n-1\}$ and $\varepsilon \in \{\pm 1\}$. The normal form of 
$\rho(\sigma_k^\varepsilon)(u)$ with respect to the decomposition $G=H_1 \ast G'$ also starts with 
$h_1^{-1}$ and ends with $h_1$.
\end{lemma}

\Proof
Let $v\in H_1 \ast H_2$. Suppose that the normal form of $v$ is
$$
v=\phi_1(x_1)\, \phi_2(y_1)\, \dots\, \phi_1(x_l)\, \phi_2(y_l),
$$
where $x_1, \dots, x_l, y_1, \dots, y_{l-1} \in H\setminus \{1\}$, and $y_l \in H$. Then
$$
\rho(\sigma_1)(v)=h_1^{-1} \cdot \phi_2(x_1) \cdot h_1^2 \phi_1(y_1) h_1^{-2} \cdot \dots \cdot 
\phi_2(x_l) \cdot h_1^2 \phi_1(y_l)\ h_1^{-1},
$$
thus the normal form of $\rho(\sigma_1)(v)$ starts with $h_1^{-1}$.

Similarly, if the normal form of $v$ is
$$
v=\phi_2(y_1)\, \phi_1(x_1)\, \dots\, \phi_2(y_l)\, \phi_1(x_l),
$$
where $x_1, \dots, x_l, y_2, \dots, y_l \in H \setminus \{1\}$ and $y_1 \in H$,
then the normal form of $\rho(\sigma_1)(v)$ ends with $h_1$.

Now, write
$$
u=v_0\, w_1\, v_1 \dots w_l\, v_l
$$
where $v_i \in (H_1 \ast H_2) \setminus \{1\}$ and $w_j\in (H_3 \ast \dots \ast H_n) \setminus 
\{1\}$, and $l\ge 0$. The hypothesis that $u$ starts with $h_1^{-1}$ implies that $v_0$ starts 
with $h_1^{-1}$, and the hypothesis that $u$ ends with $h_1$ implies that $v_l$ ends with $h_1$. 
Both groups, $H_1 \ast H_2$ and $H_3 \ast \dots \ast H_n$, are invariant by $\rho(\sigma_1)$,
and $\rho(\sigma_1)$ is the identity on $H_3\ast \dots \ast H_n$. So,
$$
\rho(\sigma_1)(u)= \rho(\sigma_1)(v_0) \cdot w_1\cdot \rho(\sigma_1)(v_1) \cdot \dots w_l \cdot 
\rho(\sigma_1)(v_l).
$$
By the above observations, $\rho(\sigma_1)(v_0)$ starts with $h_1^{-1}$ and $\rho(\sigma_1)(v_l)$ 
ends with $h_1$, thus $\rho(\sigma_1)(u)$ starts with $h_1^{-1}$ and ends with $h_1$.

Let $k \in \{2, \dots, n-1\}$ and $\varepsilon\in \{\pm 1\}$. Write
$$
u=h_1^{-1}\, w_1\, v_1 \dots v_{l-1}\, w_l\, h_1,
$$
where $v_1, \dots, v_{l-1} \in H_1 \setminus \{1\}$ and $w_1, \dots, w_l \in G'\setminus \{1\}$. 
Both groups, $H_1$ and $G'$, are invariant by $\rho(\sigma_k^\varepsilon)$, and 
$\rho(\sigma_k^\varepsilon)$ is the identity on $H_1$. So,
$$
\rho(\sigma_k^\varepsilon)(u)= h_1^{-1}\cdot \rho(\sigma_k^\varepsilon)(w_1) \cdot v_1 \cdot 
\dots \cdot v_{l-1} \cdot \rho(\sigma_k^\varepsilon)(w_l) \cdot h_1,
$$
thus the normal form of $\rho(\sigma_k^\varepsilon)(u)$ starts with $h_1^{-1}$ and ends with 
$h_1$.
\endproof

\paragraph{Proof of Proposition \ref{4.1}.} 
We argue by induction on $n$. Assume $n=2$. We have
$$\displaylines{
\rho(\sigma_1^{2l})(h_1)=(h_2h_1)^{-l} h_1 (h_2h_1)^l \neq h_1, \quad {\rm for}\ l \in \Z
\setminus \{0\} \cr
\rho(\sigma_1^{2l+1})(h_1)=(h_2h_1)^{-l} h_1^{-1}h_2h_1 (h_2h_1)^l \neq h_1, \quad {\rm for}\ 
l \in \Z\cr}
$$
thus the representation $\rho: B_2 \to \Aut(H_1 \ast H_2)$ is faithful.

Now, assume $n\ge 3$. Let $\beta \in B_n\setminus \{1\}$. By Proposition \ref{4.2}, either $\beta \in 
B_{n-1}$, or $\beta=\alpha_0 \sigma_1 \dots \sigma_1 \alpha_l$, where $l\ge 1$ and $\alpha_0, 
\dots, \alpha_l \in B_{n-1}$, or $\beta=\alpha_0 \sigma_1^{-1} \dots \sigma_1^{-1} \alpha_l$, 
where $l\ge 1$ and $\alpha_0, \dots, \alpha_l \in B_{n-1}$.

Suppose $\beta \in B_{n-1}$. By induction, $\rho(\beta)$ acts non-trivially on $G'=H_2 \ast \dots 
\ast H_n$, thus $\rho(\beta)$ acts non-trivially on $G=H_1 \ast G'$.

Suppose $\beta=\alpha_0 \sigma_1 \dots \sigma_1 \alpha_l$, where $l\ge 1$ and $\alpha_0, 
\dots, \alpha_l \in B_{n-1}$. Let
$$
u=\rho(\sigma_1 \alpha_l)(h_1)=\rho(\sigma_1)(h_1)= h_1^{-1} h_2 h_1.
$$
By Lemma \ref{4.3}, the normal form of $\rho(\alpha_0 \sigma_1 \dots \sigma_1 \alpha_{l-
1})(u)=\rho(\beta)(h_1)$ starts with $h_1^{-1}$ and ends with $h_1$. In particular, 
$\rho(\beta)(h_1)\neq h_1$, thus $\rho(\beta)\neq \Id$.

Suppose $\beta=\alpha_0 \sigma_1^{-1} \dots \sigma_1^{-1} \alpha_l$, where $l\ge 1$ and 
$\alpha_0, \dots, \alpha_l \in B_{n-1}$. By the previous case, $\rho(\beta^{-1}) \neq \Id$, thus 
$\rho(\beta) \neq \Id$.
\endproof


\section{Garside groups}\label{Sect5}

Our objectives in this section are twofold. 
Firstly we give a brief presentation of the definition and
salient properties of a Garside group. 
Secondly we establish a criterion (or set of criteria)
which allows one to show that a group given by a certain type of presentation is indeed
a Garside group, and which we will make use of in the subsequent section.
Our presentation of the subject draws in many ways from the work of Dehornoy
\cite{Deh3, Deh5} as well as \cite{DP}, and, like all treatments of Garside groups,
is inspired ultimately by the seminal papers of Garside \cite{Gar}, 
on braid groups, and Brieskorn and Saito \cite{BS}, on Artin groups. 
 
\begin{defn}
Let $M$ be an arbitrary monoid. 
We say that $M$ is {\it atomic} if there exists a function $\nu: M \to \N$ such that
\begin{description}
\item{$\bullet$} $\nu(a)=0$ if and only if $a=1$;
\item{$\bullet$} $\nu(ab) \geq \nu(a) + \nu(b)$ for all $a,b \in M$.
\end{description}
Such a function $\nu: M \to \N$ is called a {\it norm} on $M$.

An element $a\in M$ is called an \emph{atom} if it is indecomposable, namely, if $a=bc$ then
either $b=1$ or $c=1$.
\end{defn}

This definition of atomicity  is taken from \cite{DP}.  See 
\cite{DP}, Proposition 2.1, for a list of further properties all equivalent to atomicity. 
In the same paper it is shown that any generating set of $M$ contains the set of all atoms.
In particular, $M$ is finitely generated if and only if it has only finitely many atoms. 

Given that a monoid $M$ is atomic, we may define left and right invariant 
partial orders $\leq_L$ and $\leq_R$ on $M$ as follows:
\begin{description}
\item{$\bullet$} set $a\leq_Lb$ if there exists $c \in M$ such that $ac=b$;
\item{$\bullet$} set $a\leq_Rb$ if there exists $c\in M$ such that $ca=b$.
\end{description}
We shall call these the \emph{left} and \emph{right divisibility orders} on $M$.

\begin{defn}\label{GarM}
A {\it Garside monoid} is a monoid $M$ such that
\begin{description}
\item{(i)} $M$ is atomic and finitely generated;
\item{(ii)} $M$ is cancellative;
\item{(iii)} $(M,\leq_L)$ and $(M,\leq_R)$ are lattices;
\item{(iv)} there exists an element $\De\in M$, which we call a \emph{Garside element}, such that
	\begin{description}
	\item{(a)} the set $L(\De):=\{x\in M : x\le_L\De\}$ generates $M$, and
	\item{(b)} the sets $L(\De)$ and $R(\De):= \{x\in M : x\le_R\De\}$ are equal.
	\end{description}
\end{description}
\end{defn}

\Remark 
Elsewhere in the literature 
the condition that $M$ is finitely generated is often incorporated into
condition (iv) of the definition by saying that the set $L(\De)$ is finite.
It seems more natural to state this condition separately. Note that, if $M$
is finitely generated and atomic, then $L(a)= \{x \in M: x\le_La\}$ is finite for 
all $a \in M$.

\begin{defn}
For any monoid $M$ one can define the group $G(M)$ which is 
presented by the generating set $M$ and relations 
$ab=c$ whenever $ab=c$ in $M$. There is an obvious  canonical homomorphism $M\to G(M)$. 
This homomorphism is not injective in general.
The group $G(M)$ is known as the {\it group of fractions} of $M$. Define a {\it Garside group} to 
be the group of fractions of a Garside monoid.
\end{defn}

\Remark
\begin{description}
\item{(1)} A Garside monoid $M$ satisfies \"Ore's conditions, thus the canonical homomorphism $M 
\to G(M)$ is injective. Moreover the partial order $\leq_L$ (resp. $\leq_R$) extends to a 
left invariant (resp. right invariant) lattice order on $G(M)$ with positive cone $M$.
\item{(2)} A Garside element is never unique. For example, if $\Delta$ is a Garside element, then 
$\Delta^k$ is also a Garside element for all $k \ge 1$ (see \cite{Deh5}, Lemma 2.2).
\end{description}

Let $M$ be a Garside monoid. The lattice operations of $(M,\le_L)$ are denoted by $\vee_L$ and 
$\wedge_L$. For $a,b \in M$, we denote by $a\backslash_L b$ the unique element of $M$ such that 
$a (a\backslash_Lb)= a\vee_L b$. Similarly, the lattice operations of $(M,\le_R)$ are denoted by 
$\vee_R$ and $\wedge_R$, and, for $a,b \in M$, we denote by $b/_Ra$ the unique element of $M$ such 
that $(b/_Ra)a=a\vee_R b$.

\medskip
Now, before establishing our criterion for a group to be a Garside group, we briefly explain how 
to define a biautomatic structure on a given Garside group. By \cite{EpAl}, such a structure 
furnishes solutions to the word problem and to the conjugacy problem, and it implies that the 
group has quadratic isoperimetric inequalities. We refer to \cite{EpAl} for definitions and 
properties of automatic groups, and to \cite{Deh5} for more details on the biautomatic structures 
on Garside groups.

Let $M$ be a Garside monoid, and let $\Delta$ be a Garside element of $M$. For $a \in M$, we 
write $\pi_L(a)= \Delta \wedge_L a$ and denote by $\partial_L(a)$ the unique element of $M$ such 
that $a=\pi_L(a) \, \partial_L(a)$. Using the fact that $M$ is atomic and that $L(\Delta)=\{x\in 
M : x\le_L \Delta\}$ contains all the atoms, one can easily show that $\pi_L(a) \neq 1$ if $a 
\neq 1$, and that there exists some positive integer $k$ such that $\partial_L^k(a)=1$. Let $k$ be the 
lowest integer satisfying $\partial_L^k(a)=1$. Then the expression
$$
a= \pi_L(a)\cdot \pi_L (\partial_L (a))\cdot \dots\cdot \pi_L( \partial_L^{k-1}(a))
$$
is called the {\it normal form} of $a$.

Let $G=G(M)$ be the group of fractions of $M$. Let $c\in G$. Since $G$ is a lattice
with positive cone $M$ the element $c$ can be written $c=a^{-1}b$ with $a,b \in M$.
Obviously, $a$ and $b$ can be chosen so that $a \wedge_Lb=1$ and, with this extra
condition, are unique. Let $a=a_1a_2 \dots a_p$ and $b=b_1 b_2 \dots b_q$ be the
normal forms of $a$ and $b$, respectively. Then the expression
\[
c= a_p^{-1} \dots a_2^{-1} a_1^{-1} b_1 b_2 \dots b_q
\]
is called the {\it normal form} of $c$. 

The following result can be found in \cite{Deh5}, Section 3.

\begin{thm}[Dehornoy]\label{new5.1}
Let $M$ be a Garside monoid and let $G$ be the group of fractions of $M$. Then the normal forms 
of the elements of $G$ form a symmetric rational language on the (finite) set $L(\Delta)$ which
has the fellow traveler property. In particular, $G$ is biautomatic.\noproof
\end{thm}

We turn now to establish our criterion.

For a finite set $S$, we denote by $S^\ast$ the free monoid on $S$. The elements of $S^\ast$ are 
called {\it words} on $S$. The empty word is denoted by $\epsilon$. Let $\equiv$ be a congruence 
relation on $S^\ast$, and let $M= (S^\ast/\equiv)$. For $w\in S^\ast$, we denote by $\ov{w}$ the 
element of $M$ represented by $w$, and we call $w$ an {\it expression} of $\ov{w}$.

\begin{defn}
A {\it complement}  is a function $f:S \times S\to S^\ast$ such that $f(x,x)=\epsilon$
for all $x\in S$. To a complement $f: S \times S \to S^\ast$ we associate the following two monoids.
$$
\begin{array}{rl}
M_L^f &= \langle S\ |\ xf(x,y) = yf(y,x)\ {\rm for}\ x,y\in S \rangle^+,\\
M_R^f &= \langle S\ |\ f(y,x)x = f(x,y)y\ {\rm for}\ x,y\in S \rangle^+.\\
\end{array}
$$
For $u,v\in S^\ast$, we write $u\equiv_L^f v$ if $u$ and $v$ are expressions of the same element 
of $M_L^f$, and we write $u\equiv_R^f v$ if $u$ and $v$ are expressions of the same element of 
$M_R^f$.
\end{defn}

\begin{defn}
A word $w$ in $(S \cup S^{-1})^\ast$ is {\it $f$-reversible on the left in one 
step} to a word $w'$ if $w'$ is obtained from $w$ by replacing some subword $x^{-1}y$ (with 
$x,y\in S$) by the corresponding word $f(x,y)f(y,x)^{-1}$. Let $p\ge 0$. We say that $w$ is {\it 
$f$-reversible on the left in $p$ steps} to a word $w'$ if there exists a sequence $w_0=w, w_1, 
\dots, w_p=w'$ in $(S \cup S^{-1})^\ast$ such that $w_{i-1}$ is $f$-reversible on the left in one 
step to $w_i$ for all $i=1, \dots, p$. The property ``$w$ is $f$-reversible on the left to 
$w'$'' is denoted by $w\rightarrow_L^f w'$.

We define the {\it $f$-reversibility on the right} in a similar way, replacing subwords $yx^{-
1}$ (with $x,y\in S$) by the corresponding words $f(x,y)^{-1} f(y,x)$. The property ``$w$ is 
$f$-reversible on the right to $w'$'' is denoted by $w \rightarrow_R^f w'$.
\end{defn}

\bigskip
It is shown in [Deh3] that a reversing process is confluent, namely:

\begin{prop}[Dehornoy, \cite{Deh3}, Lemma 1]\label{New5.2} 
 Let $f: S \times S \to S^\ast$ be a 
complement, and let $w \in (S \cup S^{-1})^\ast$. Suppose that the word $w$ is $f$-reversible on 
the left in $p$ steps to a word $uv^{-1}$, with $u,v \in S^\ast$. Then any sequence of left $f$-reversing 
transformations starting from $w$ leads in $p$ steps to $uv^{-1}$.
\noproof
\end{prop}

\begin{defn}
Let $f:S \times S \to S^\ast$ be a complement and let $u,v \in S^\ast$. Assume 
that there exist $u',v' \in S^\ast$ such that $u^{-1}v \rightarrow_L^f u' (v')^{-1}$. By 
Proposition \ref{New5.2}, $u'$ and $v'$ are unique (if they exist). Then we write $u'=C_L^f(u,v)$ and 
$v'=C_L^f(v,u)$. One has
$$
u C_L^f(u,v) \equiv_L^f v C_L^f (v,u)
$$
(see [Deh3, Lem. 2]). If no such words $u',v'$ exist then we write  $C_L^f(u,v)=C_L^f(v,u)=\infty$.

Similarly, define the words $C_R^f(u,v)$ and $C_R^f(v,u)$ to be the unique elements of 
$S^\ast$ which satisfy $vu^{-1} \rightarrow_R^f C_R^f(u,v)^{-1} C_R^f(v,u)$, or write 
 $C_R^f(u,v)=C_R^f(v,u)=\infty$ if no such words exist.
\end{defn}

\begin{defn}[Dehornoy, \cite{Deh3}, p.120]
Let $f: S \times S \to S^\ast$ be a complement.
We say that $f$ is {\it coherent on the left} if, for all $x,y,z\in S$ such that 
$C_L^f(f(x,y), f(x,z))\neq\infty$ we have
\[
C_L^f(f(x,y), f(x,z))\equiv_L^f C_L^f(f(y,x), f(y,z))\,.
\]
Similarly, we say that $f$ is {\it coherent on the right} if, for all $x,y,z\in S$ such that 
$C_R^f(f(z,x), f(y,x))\neq\infty$ we have 
\[
C_R^f(f(z,x),f(y,x))\equiv_R^f C_R^f(f(z,y),f(x,y))\,.
\]
\end{defn}

\medskip
A partially ordered set $(X,\leq)$ is said to be a \emph{quasi-lattice}, or \emph{quasi-lattice ordered},
if every pair of elements $x,y\in X$ which has a common upper bound ($z$ such that $x\leq z$ and $y\leq z$)
has a least upper bound, usually written $x\vee y$. 

The proof of the following proposition can be more or less  reconstructed from 
Garside's original treatment of the braid monoids \cite{Gar}, or the similar
treatment of Artin monoids in
\cite{BS}. The result appears in almost precisely this form (with some notational differences)
as Lemma 4 of \cite{Deh3}.  

\begin{prop}\label{qlo}
Let $M$ be an atomic monoid with generating set $S$,
and suppose that $f: S \times S \to S^\ast$ is a
complement which is coherent on the left and such that $M=M_L^f$. 
Then the following holds:
\begin{description}
\item[(LCQL)] For all $u,v,x,y\in S^\ast$ such that $ux\equiv_L^f vy$,
there exists $w\in S^\ast$ such that $x\equiv_L^f C_L^f(u,v)w$ and $y\equiv_L^f C_L^f(v,u)w$.
\end{description}
In particular, {\bf (LCQL)} implies that $M$ is left cancellative and $(M,\leq_L)$
is a quasi-lattice. 
\end{prop}

\Proof
We refer the reader to \cite{Deh3}, Lemma 4, for the proof of the statement {\bf (LCQL)}.
The fact that $M$ is left cancellative comes from putting $u=v$ in {\bf (LCQL)}.
Also, if $u$ and $v$ represent elements $\ov u$ and $\ov v$ respectively, and if $\ov u$ and $\ov v$
have a common upper bound (represented by words $ux\equiv vy$ for some $x,y\in S^\ast$), then the
least upper bound $\ov u \vee_L \ov v$ is the element represented by $uC_L^f(u,v)\equiv vC_L^f(v,u)$. 
The statement {\bf (LCQL)} implies that this element divides all common upper bounds of $\ov u$ and $\ov v$.
\endproof

Now, Proposition \ref{qlo}, together with \cite{DP} and \cite{Deh5}, 
permit the following criterion for a monoid $M$ to be a Garside monoid:

\begin{crit}\label{criterion}
Let $M$ be a monoid. Then $M$ is a Garside monoid if and only if it satisfies the following properties:

\begin{description}
\item{(C1)} $M$ is finitely generated and atomic; 
\item{(C2)} there exist complements $f: S_1 \times S_1 \to S_1^\ast$, coherent on the left,
 and  $g: S_2 \times S_2 \to S_2^\ast$, coherent on the right, such that $M\cong M_L^f$ and $M\cong M_R^g$;
\item{(C3)} $M$ possesses a Garside element, namely an element $\De\in M$ such that every atom 
of $M$ left divides $\De$ and 
the sets $L(\Delta) = \{x \in M: x\le_L \Delta\}$ and $R(\Delta)= \{ x \in M: x\le_R \Delta\}$
are equal.
\end{description}
\end{crit}

\Proof
Let $M$ be a Garside monoid. Clearly, $M$ satisfies (C1) and (C3). So, we just need to show that 
$M$ satisfies (C2). Choose some finite generating set $S$ for $M$, and consider complements 
$f: S \times S \to S^\ast$ and $g: S \times S \to S^\ast$ such that
$$
\ov{f(x,y)} = x\backslash_L y\,,\quad \ov{g(x,y)}= y/_Rx\,,
$$
for all $x,y\in S$. Then, by \cite{DP}, Theorem 4.1, one has $M=M_L^f=M_R^g$, and, by 
\cite{Deh5}, Lemma 5.2, $f$ is coherent on the left and $g$ is coherent on the right.

Now, recall the statement of \cite{Deh5}, Proposition 2.1.

Suppose that $M$ is a monoid which satisfies the following properties:
\begin{description}
\item{(D1)} $M$ is finitely generated and atomic;
\item{(D2)} $M$ is left and right cancellative;
\item{(D3)} $(M, \le_L)$ is a quasi-lattice;
\item{(D4)} there exists a finite subset $P \subset M$ which generates $M$ and which is closed 
under $\backslash_L$ (namely, if $a,b \in P$, then $a \backslash_L b \in P$).
\end{description}
Then $M$ is a Garside monoid.

Let $M$ be a monoid which satisfies (C1), (C2), (C3). 
We wish to show that $M$ is a Garside monoid.
By Proposition \ref{qlo}, $M$ satisfies 
(D1), (D2) and (D3). So, it remains to show that $M$ satisfies (D4). Let $P=L(\Delta) = 
R(\Delta)$. Note that, by hypothesis, $P$ generates $M$. Let $a,b \in P$. Since $a \le_L \Delta$ 
and $b \le_L \Delta$, we have $a \vee_L b \le_L \Delta$. Let $c \in M$ such that $\Delta = (a 
\vee_L b)c= a (a\backslash_L b) c$. Then $(a\backslash_L b) c \le_R \Delta$, thus $(a\backslash_L 
b) c \le_L \Delta$ (since $L(\Delta) = R(\Delta)$), therefore $(a\backslash_L b) \le_L \Delta$, 
that is $(a\backslash_L b) \in P$.
\endproof

\Remark
In the context of Garside groups, the reversing processes are used not only to determine whether 
a group is a Garside group, but it is also a very useful tool for solving the word problem and 
to explicitly compute normal forms. For instance, if $M$ is a Garside monoid, then one can find 
complements $f: S \times S \to S^\ast$ and $g: S \times S \to S^\ast$ such that $M=M_L^f 
=M_R^g$. If $w$ is in $(S \cup S^{-1})^\ast$, then any sequence of right $g$-reversing 
transformations leads to a word $u^{-1} v$ where $u, v \in S^\ast$, and one has $\ov{w} = 
\ov{u}^{-1}\, \ov{v}$. On the other hand, if $u,v \in S^\ast$, then 
$\ov{u} \backslash_L \ov{v}$ is represented by $C_L^f(u,v)$,
$\ov{u} \vee_L \ov{v}$ is 
represented by $u C_L^f(u,v)$, and $\ov{u} \wedge_L \ov{v}$ can be computed by means of the 
equality  
$$
\ov{u} \wedge_L \ov{v} = (\ov{u} \vee_L \ov{v}) /_R ((\ov{u} \backslash_L \ov{v}) \vee_R (\ov{v} 
\backslash_L \ov{u}))
$$
(see \cite{Deh5}, Lemma 2.6).

\medskip
In the next section we shall need the following characterization of a Garside element.

\begin{lemma}[Garside elements]\label{Deltalemma}
 Let $M$ be a cancellative atomic monoid with atom set $A$, and suppose that 
$\De\in M$ is such that $L(\De):=\{x\in M : x\le_L\De\}$ generates $M$. Define 
$R(\De):=\{x\in M : x\le_R\De\}$. Then the following are equivalent:
\begin{description}
\item{(1)} $L(\De)=R(\De)\,$;
\item{(2)} $A.\De=\De.A\,$;
\item{(3)} $M.\De=\De.M\,$;
\item{(4)} there exists a monoid automorphism $\tau :M\to M$ such that $w\De=\De\tau(w)$ for
all $w\in M$. (In particular, $\tau(A)=A$. Also $\tau$ is necessarily unique.)
\end{description}
\end{lemma}

\Proof 
By cancellativity, there is a well-defined bijection $c:L(\De)\to R(\De)$ such that
$x.c(x)=\De$ for all $x\in L(\De)$. Suppose that (1) holds. Then $c$ is a bijection $L(\De)\to L(\De)$ and
we may define $\tau=c^2$, also a bijection of $L(\De)\to L(\De)$. Note that for $x\in L(\De)$,
we also have $c(x)\in L(\De)$, so that $\De$ may be written $c(x)c^2(x)$.
Therefore, for all $x\in L(\De)$, we have $x.\De=x.c(x).c^2(x)=\De.\tau(x)$.
Since $L(\De)$ generates $M$, it follows that $\De\leq_L w\De$ for all $w\in M$ (in fact,
if $w=x_1x_2...x_n$ with $x_i\in L(\De)$, then $w\De=\De\tau(x_1)...\tau(x_n)$). 
By left cancellativity, there is therefore a unique well-defined function $\tau: M\to M$ such that
$w\De=\De\tau(w)$ for all $w\in M$. By right cancellativity, $\tau$ must be injective.
Moreover, given $x,y,z\in M$ such that $z=xy$, we have $\De\tau(z)=z\De=x\De\tau(y)=\De\tau(x)\tau(y)$
and, by cancellation, $\tau(z)=\tau(x)\tau(y)$. Thus $\tau$ is a monoid homomorphism.
By a similar argument, we may construct the inverse homomorphism $\tau^{-1}$
in order to show that $\tau$  is in fact an automorphism of $M$. Thus (1) implies (4).

Now suppose that (2) holds: $A\De=\De A$. Then by left and right cancellativity, there is 
a well-defined bijection $\tau:A\to A$ such that $a\De=\De\tau(a)$ for all $a\in A$. As in the previous
paragraph this extends to an automorphism of $M$ such that $w\De=\De\tau(w)$ for all
$w\in M$. Thus (2) implies (4). By the same reasoning one can show that (3) implies (4).
On the other hand, both (2) and (3) are obvious consequences of (4).

Finally, we show that (4) implies (1). Suppose that (4) holds.
In particular, we have $\tau(\De)=\De$. Therefore,  $x\leq_L\De$ if and only if $\tau(x)\leq_L\Delta$
(since $\tau$ is a monoid automorphism). In other words, $\tau(L(\De))=L(\De)$.
On the other hand, the equation $x\De=\De\tau(x)$ shows, by left cancellation, that 
if $x\in L(\De)$ then $\tau(x)\in R(\De)$, and, by right cancellation, that
if  $y=\tau(x)\in R(\De)$ then $\tau^{-1}(y)\in L(\De)$. Thus $\tau(L(\De))=R(\De)$. 
But then $\tau(L(\De))=L(\De)=R(\De)$, giving (1).
\endproof


\section{Semi-direct products}\label{Sect6}

We turn back to the Artin type representations. Let $H$ be a group, let $h \in H \setminus \{1\}$, 
let $G=H_1 \ast \dots \ast H_n$, where $H_i$ is a copy of $H$, and let $\rho: B_n \to \Aut(G)$ be 
the Artin type representation associated to $(H,h)$. The aim of this section is to prove the 
following.

\begin{thm}\label{6.1}
Assume that $H$ is the group of fractions of a Garside monoid $M$ and that $h$ is a Garside 
element. Let $\wtil{G} = G \rtimes_\rho B_n$, and let $\wtil M$ be the submonoid of $\wtil G$ 
generated by $M_1=\phi_1(M)$ and the monoid $B_n^+$ of positive braids. Then $\wtil M$ is a 
Garside monoid, $\Delta=( h_1 \sigma_1 \sigma_2 \dots \sigma_{n-1})^n$ is a Garside element of 
$\wtil M$, and $\wtil G$ is the group of fractions of $\wtil M$.
\end{thm}

\bigskip
The first step in the proof of Theorem 6.1 consists on finding a presentation for $G 
\rtimes_\rho B_n$, namely:

\begin{prop}\label{presentation}
Let $H=\langle S \mid \Cal R\rangle$ be a presentation for $H$, and let $D 
\in S^\ast$ be an expression of $h$. Then $\wtil G= G \rtimes_\rho B_n$ has a 
presentation with generators
\[
S \cup \{\sigma_1, \dots, \sigma_{n-1}\},
\]
and with relations
\[
\begin{array}{cl}
r&\text{for }\ r \in \Cal R\,,\\
\sigma_i \sigma_{i+1} \sigma_i = \sigma_{i+1} \sigma_i \sigma_{i+1}&\text{for }\ i=1, \ldots, n-2\,,\\
\sigma_i \sigma_j = \sigma_j \sigma_i &\text{for }\ |i-j| \ge 2\,,\\
\sigma_i x = x \sigma_i &\text{for }\ x \in S \text{ and } i=2, \ldots , n-1\,,\\
x \sigma_1 D \sigma_1 = \sigma_1 D \sigma_1 D^{-1} x D &\text{for }\ x\in S\,.
\end{array}
\]
\end{prop}

\Proof
Let $\wtil G_0$ denote the abstract group generated by $S \cup \{ \sigma_1, \dots, 
\sigma_{n-1} \}$ and subject to the relations given in the statement of Proposition \ref{presentation}.
Let $X=(\cup_{i=1}^n \phi_i(S)) \cup \{\sigma_1, \dots, \sigma_{n-1}\}$. With a little effort one can
verify that the mapping $\varphi: X \to \wtil G_0$ defined by
\[
\begin{array}{rll}
\varphi(\phi_i(x))&=\ \sigma_{i-1}^{-1} \dots \sigma_1^{-1} D^{i-1} x D^{1-i}\sigma_1 \dots \sigma_{i-1}
&\text{for }\ i=1, \dots, n \text{ and } x \in S\\
\varphi(\sigma_i)&=\ \sigma_i &\text{for }\ i=1, \dots, n-1\\
\end{array}
\]
determines a homomorphism $\varphi: \wtil G \to \wtil G_0$, and somewhat more easily 
that the mapping $\psi: S \cup \{ \sigma_1, \dots, \sigma_{n-1} \} \to \wtil G$ defined by
\[
\begin{array}{rll}
\psi(x) &=\ \phi_1(x) &\text{for }\ x\in S\\
\psi(\sigma_i) &=\ \sigma_i &\text{for }\ i=1,\ldots ,n-1\\
\end{array}
\]
determines a homomorphism $\psi: \wtil G_0 \to \wtil G$. One checks without too much difficulty
that $(\psi \circ \varphi)(a)=a$ for all $a \in X$, and $(\varphi \circ \psi)(b) =b$ for
all $b \in S \cup \{ \sigma_1, \dots, \sigma_{n-1}\}$, thus $\psi \circ \varphi = 
\Id_{\wtil G}$ and $\varphi \circ \psi = \Id_{\wtil G_0}$.
\endproof

\paragraph
{Proof of Theorem \ref{6.1}.}
Let $\tau:M\to M$ denote the
automorphism of $M$ induced by conjugation by $h^{-1}$, so that $xh=h\tau(x)$ for all $x\in M$ 
(see Lemma \ref{Deltalemma}). Let $S$ be a finite generating set 
for $M$. We may, and do, choose $S$ so that $\tau(S)=S$
(for instance we may simply choose $S$ to be the set of atoms of $M$).
Define $f: S \times S \to S^\ast$ such that $\ov{xf(x,y)}=\ov{yf(y,x)}=x\vee_L y$ for all pairs
$x,y\in S$. Similarly define $g: S \times S \to S^\ast$ such that 
$\ov{g(x,y)y}=\ov{g(y,x)x}=x\vee_R y$ for  all pairs $x,y\in S$. 
As pointed out in the proof of Criterion \ref{criterion}, one has $M=M_L^f=M_R^g$,
$f$ is coherent on the left, and $g$ is coherent on the right.
We simply write $\sim$
for the congruence relation on $S^\ast$ defined by the relations in $M$ (namely, $\equiv_L^f$, or equally
$\equiv_R^g$).  Let $D \in S^\ast$ be an expression of $h$. 
Note that for $x\in S$ we have $xD\sim D\tau(x)$ and $\tau^{-1}(x)D\sim Dx$, where $\tau(x)$
and $\tau^{-1}(x)$ also denote elements of the generating set $S$.
The last family of relations appearing in Proposition \ref{presentation} may be  replaced
with  $x\sigma_1 D\sigma_1=\sigma_1 D\sigma_1\tau(x)$ for all $x\in S$, or equivalently with
$\tau^{-1}(x)\sigma_1 D\sigma_1=\sigma_1 D\sigma_1 x$ for all $x\in S$.

Let $X=S \cup \{ \sigma_1, \dots, \sigma_{n-1}\}$. Let $F: X \times X \to X^\ast$ be the 
complement defined by
\[
\begin{array}{rllrll}
F(x,y) &=\ f(x,y) &\text{for }\ x,y\in S 
&F(\sigma_i,x) &=\ x &\text{for }\ x\in S\text{ and }i\geq 2\\
F(x,\sigma_1) &=\ \sigma_1 D \sigma_1 &\text{for }\ x \in S
&F(\sigma_i, \sigma_j) &=\ \sigma_j \sigma_i &\text{for }\ |i-j|=1\\
F(\sigma_1,x) &=\ D \sigma_1 \tau(x) &\text{for }\ x \in S
&F(\sigma_i, \sigma_j) &=\ \sigma_j &\text{for }\ |i-j| \ge 2\\
F(x,\sigma_i) &=\ \sigma_i &\text{for }\ x \in S\text{ and }i\geq 2 &&&
\end{array}
\]
and let $G: X \times X \to X^\ast$ be the complement defined by
\[
\begin{array}{rllrll}
G(x,y) &=\ g(x,y) &\text{for }\ x,y\in S 
&G(x,\sigma_i) &=\ x &\text{for }\ x\in S\text{ and }i\geq 2\\
G(\sigma_1,x) &=\ \sigma_1 D \sigma_1 &\text{for }\ x \in S
&G(\sigma_j, \sigma_i) &=\ \sigma_i \sigma_j &\text{for }\ |i-j|=1\\
G(x,\sigma_1) &=\ \tau^{-1}(x) \sigma_1 D &\text{for }\ x \in S
&G(\sigma_j, \sigma_i) &=\ \sigma_j &\text{for }\ |i-j| \ge 2\\
G(\sigma_i,x) &=\ \sigma_i &\text{for }\ x \in S\text{ and }i\geq 2 &&&
\end{array}
\]
Let $\wtil M_0$ denote the monoid defined by the presentation with generators $X$ and relations
as laid out in Proposition \ref{presentation}. Then clearly $\wtil M_0\cong M_L^F\cong M_R^G$.
We denote by $\approx$ the congruence relation on $X^\ast$ defined by the relations of $\wtil M_0$.
(So $\approx$ is the same congruence relation as $\equiv_L^F$ and $\equiv_R^G$).
We proceed now to show that $\wtil M_0$ satisfies the Criterion \ref{criterion}
with complements $F$ and $G$ and Garside element $\De=(D\sig_1\sig_2\ldots\sig_{n-1})^n$.
It follows that $\wtil M_0$ is a Garside monoid with group of fractions $\wtil G$
and is canonically isomorphic to the submonoid $\wtil M\subset \wtil G$ in the statement
of the Theorem.

Clearly $\wtil M_0$ is finitely generated. We check that $\wtil M_0$ is atomic.
Let $\nu:M\to\N$ be a norm for $M$.  Let $\Sigma = \{ \sigma_1, \dots, \sigma_{n-1}\}$ and 
define the function $\ell: \Sigma^\ast \to \N$  by $\ell(\sigma_{i_1} \dots \sigma_{i_l}) = l$.
We define a function $\wtil \nu: X^\ast \to \N$ as follows. Let $w \in X^\ast$. 
Write $w=u_1v_1 \dots u_lv_l$, where $u_1 \in 
S^\ast$, $u_2, \dots, u_l \in S^\ast\setminus \{ \epsilon\}$, $v_1, \dots, v_{l-1} \in 
\Sigma^\ast \setminus \{ \epsilon \}$, and $v_l \in \Sigma^\ast$.
Then
$$
\wtil \nu (w) = \nu(u_1u_2...u_l) + \ell(v_1 v_2 \dots v_l).
$$
One can easily verify that $\wtil\nu$ is invariant with respect to all of the relations
given in Proposition \ref{presentation}, and therefore defines a function $\wtil\nu :\wtil M_0\to \N$.
Moreover, it is easily seen that $\wtil\nu$ is a norm, and therefore $\wtil M_0$ is atomic.

\begin{figure}[ht]
\begin{center}
\includegraphics[width=14cm]{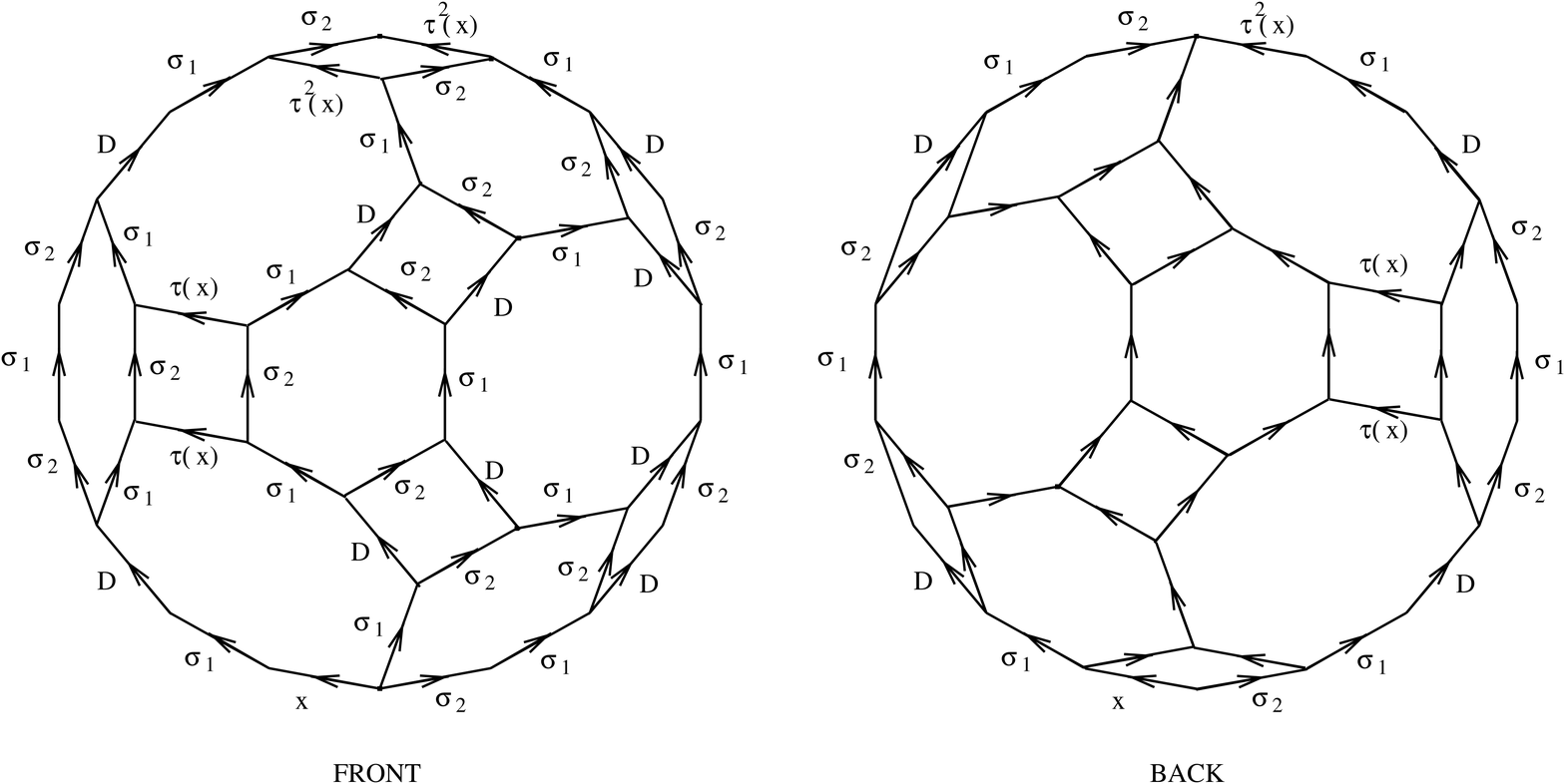}
\end{center}
\caption{Left coherence of $F$ with respect to triple $\{\sig_1,\sig_2,x\}$.}\label{Pic5}
\end{figure}

\begin{figure}[ht]
\begin{center}
\includegraphics[width=6cm]{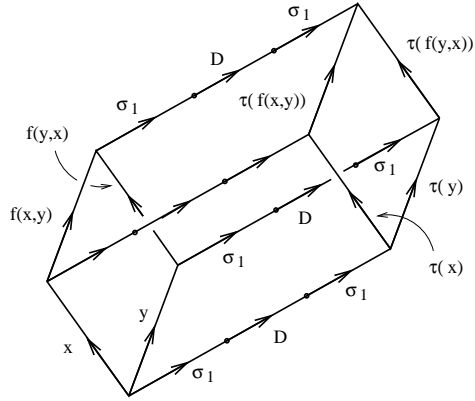}
\end{center}
\caption{Left coherence of $F$ with respect to triple $\{\sig_1,x,y\}$.}\label{Pic6}
\end{figure}

The proof that $F$ is coherent on the left may be deduced from the existence, for each triple
$\al,\be,\ga\in X$, of a certain tiling of the 2-sphere by relations from $M_L^F$
(i.e: relations of the form $\al F(\al,\be)\approx\be F(\be,\al)$ for $\al,\be\in X$)
as illustrated in Figures \ref{Pic5} and \ref{Pic6}. We illustrate the two most difficult cases, namely when 
$\{\al,\be,\ga\}=\{\sig_1,\sig_2, x\}$ for some $x\in S$ (Figure \ref{Pic5}), 
and when $\{\al,\be,\ga\}=\{\sig_1, x, y\}$ for some $x\neq y \in S$ (Figure \ref{Pic6}).
In the latter case note that, if $f(x,y)$ is written $a_1a_2..a_k$ as a product of generators $a_i\in S$
then $\tau(f(x,y))\approx\tau(a_1)\tau(a_2)..\tau(a_k)$ and the face containing $f(x,y)$ and $\tau(f(x,y))$
in the Figure \ref{Pic6} decomposes into $k$ faces corresponding to the relations 
$a_i\sigma_1 D\sigma_1\approx \sigma_1 D\sigma_1\tau(a_i)$. Similarly for $f(y,x)$.
The remaining cases are easily handled since in these cases at least one of $\al,\be,\ga$
satisfies a commuting relation (explicit in the presentation $M_L^F$) with each of the others. 

The proof that $G$ is coherent on the right is similar.

Finally we show that the word $\De=(D\sig_1\sig_2\ldots\sig_{n-1})^n$ represents
a Garside element of $\wtil M_0$. We shall employ Condition (4) of Lemma \ref{Deltalemma}.
Consider the Artin monoid presentation
\[
\begin{array}{rl}
A^+(B_n)=\<\be_1,\be_2,\ldots,\be_n\mid &\be_1\be_2\be_1\be_2=\be_2\be_1\be_2\be_1\\
&\be_i\be_{i+1}\be_i=\be_{i+1}\be_i\be_{i+1}\ \text{ for }\ 2\leq i\leq n-1\\
&\be_i\be_j=\be_j\be_i\ \text{ for }\ |i-j|\geq 2\ \>^+\,.
\end{array}
\]
This monoid $A^+(B_n)$ is well-known as the Artin monoid of type $B_n$, and has 
Garside element $\De_B=(\be_1\be_2\ldots\be_n)^n$.
Clearly there exists a monoid homomorphism $A^+(B_n)\to \wtil M_0$ such that
$\be_1\mapsto D$ and $\be_i\mapsto \sig_{i-1}$ for $i=2,3,..,n$. 
Thus any relation which is observed in $A^+(B_n)$ may be deduced in $\wtil M_0$.
In particular, the fact that $\De_B$ is a Garside element in $A^+(B_n)$ implies that 
$\De$ is left divisible by $D,\sig_1,..,\sig_{n-1}$ and hence is left divisible by every element of $X$. 
It remains to verify Condition (4) of Lemma \ref{Deltalemma}, namely that there
exists an automorphism $\wtil \tau:\wtil M_0\to\wtil M_0$ such that $w\De=\De\wtil\tau(w)$ for
all $w\in\wtil M_0$.

We already know that $\Delta_B$ is central in $A^+(B_n)$. Thus we have $\sig_i\Delta=\Delta\sig_i$
for all $i=1,2,..,n-1$. We may also check (by performing the calculation in $A^+(B_n)$) that
\[
\De\approx D\,U^{n-1}\hskip5mm\text{where }\ U:=\sig_1D\sig_1.\sig_2\sig_3\cdots \sig_{n-1}\,.
\]
Recall that $\tau$ denotes the automorphism of $M$ such that, at the level of words,
 $xD\sim D\tau(x)$ for all $x\in S^\ast$. Observe also that $xU\approx U\tau(x)$ for all $x\in S^\ast$
(or more loosely speaking, for all $x\in M$).
We now define $\wtil\tau:\wtil M_0\to\wtil M_0$ such that 
\[
\begin{array}{rll}
\wtil\tau(\sig_i)&=\sig_i              &\text{for }\ i=1,2,..,n-1\\
\wtil\tau(x)     &=\tau^n(x)\hskip3mm  &\text{for all }\ x\in M\,.
 \end{array}
\]
It is easily seen that $\wtil\tau$ is a monoid isomorphism. Moreover, for all $x\in M$,
\[
\begin{array}{rl}
x\De &\approx xDU^{n-1}\\
	&\approx D\tau(x) U^{n-1}\\
	&\approx DU^{n-1}\tau^n(x) \\
	&\approx \De\wtil\tau(x)\,,
\end{array}
\]
and $\sig_i\De\approx\De\sig_i$ for all $i=1,2,..,n-1$. Thus Condition (4) of Lemma \ref{Deltalemma}
is satisfied, and $\De$ is a Garside element.\endproof

\section{Appendix}

Throughout this section, we shall denote by $F_n$ the free group of rank $n$, and by
$x_1, \dots, x_n$ some fixed basis for $F_n$.

\bigskip\noindent
{\bf Definition.} According to Shpilrain's terminology [Shp], a {\it Wada representation of type 
(1)} is an Artin type representation associated to $(\Z,h)$, where $h$ is a non-zero integer. 
Such a representation will be denoted by $\rho_h^{(1)}: B_n \to \Aut(F_n)$. It is determined by
$$
\rho_h^{(1)} (\sigma_k) (x_i) = \begin{cases}
x_i &\text{if }\, i \neq k,k+1\\
x_k^{-h} x_{k+1} x_k^h &\text{if }\, i=k\\
x_k &\text{if }\, i=k+1
\end{cases}
$$
The {\it Wada representation of type (2)} is the representation $\rho^{(2)}: B_n \to \Aut(F_n)$ 
determined by
$$
\rho^{(2)} (\sigma_k) (x_i) = \begin{cases}
x_i &\text{if }\, i\neq k,k+1\\
x_k x_{k+1}^{-1} x_k &\text{if }\, i=k\\
x_k &\text{ if }\, i=k+1\end{cases}
$$
and the {\it Wada representation of type (3)} is the representation $\rho^{(3)}: B_n \to 
\Aut(F_n)$ determined by
$$
\rho^{(3)} (\sigma_k) (x_i) = \begin{cases}
x_i &\text{if }\, i\neq k,k+1\\
x_k^2 x_{k+1} &\text{if }\, i=k\\
x_{k+1}^{-1} x_k^{-1} x_{k+1} &\text{if }\, i=k+1
\end{cases}
$$

\noindent
{\bf Proposition A.1.} {\it (1) Let $k,l \in \Z \setminus \{ 0 \}$. Then $\rho_k^{(1)}$ and 
$\rho_l^{(1)}$ are equivalent if and only if $l=\pm k$.

\smallskip
(2) $\rho^{(2)}$ and $\rho^{(3)}$ are equivalent.

\smallskip
(3) Let $k \in \Z \setminus \{0\}$. Then $\rho^{(2)}$ and $\rho_k^{(1)}$ are not equivalent.}

\bigskip
The following lemmas A.2 and A.3 are preliminary results to the proof of Proposition A.1.

\bigskip\noindent
{\bf Lemma A.2.} {\it Consider the action of $B_n$ on $F_n$ via the representation 
$\rho_h^{(1)}$. For all $i=1, \dots, n-1$, we have
$$
F_n^{\langle \sigma_i \rangle} = \langle x_1, \dots, x_{i-1}, x_{i+1}^h x_i^h, x_{i+2}, \dots, 
x_n \rangle.
$$}

\noindent
{\bf Proof.} Write $F_n = C \ast D$, where $C= \langle x_i, x_{i+1} \rangle$ and $D= \langle x_1, 
\dots, x_{i-1}, x_{i+2}, \dots, x_n \rangle$. Both groups, $C$ and $D$, are invariant by the 
action of $\sigma_i$. Moreover, $\sigma_i$ is the identity on $D$ and acts on $C$ by $x_i \mapsto 
x_i^{-h} x_{i+1} x_i^h$, $x_{i+1} \mapsto x_i$. In particular, $F_n^{\langle \sigma_i \rangle} = 
C^{\langle \sigma_i \rangle} \ast D$.

Let $u \in C^{\langle \sigma_i \rangle}$. Write
$$
u=x_i^{n_1} x_{i+1}^{m_1} \dots x_i^{n_r} x_{i+1}^{m_r},
$$
where $r \ge 1$, $m_1, \dots, m_{r-1}, n_2, \dots, n_r \in \Z \setminus \{0\}$, and $m_r, n_1 
\in \Z$. First, suppose $n_1 \neq 0$. Then
$$
\sigma_i(u) = x_i^{-h} x_{i+1}^{n_1} x_i^{m_1} \dots x_{i+1}^{n_r} x_i^{m_r+h} =u
$$
thus
$$
-h=n_1,\ n_1=m_1,\ \dots,\ n_r=m_r,\ {\rm and}\ m_r+h=0,
$$
hence $u=(x_{i+1}^h x_i^h)^{-r}$. Now, suppose $n_1=0$. Then
$$
\sigma_i(u)= x_i^{m_1-h} x_{i+1}^{n_2} x_i^{m_2} \dots x_{i+1}^{n_r} x_i^{m_r+h},
$$
thus
$$
m_1-h=0,\ m_1=n_2,\ n_2=m_2,\ \dots,\ n_r=m_r+h,\ {\rm and}\ m_r=0,
$$
hence $u=(x_{i+1}^h x_i^h)^{r-1}$.
\qed

\bigskip\noindent
{\bf Lemma A.3.} {\it Consider the action of $B_n$ on $F_n$ via $\rho_h^{(1)}$. Then $F_n^{B_n}$ 
is the cyclic subgroup of $F_n$ generated by $x_n^h \dots x_2^h x_1^h$.}

\bigskip\noindent
{\bf Proof.} Let $u \in F_n^{B_n}$. We have $u \in F_n^{\langle \sigma_i \rangle}$ for all $i=1, 
\dots, n-1$, thus, by Lemma A.2, the reduced form of $u$ satisfies the following properties:

\smallskip
$\bullet$ all the exponents are either equal to $h$ or equal to $-h$;

\smallskip
$\bullet$ if $i \neq 1$, then $x_i^h$ is followed by $x_{i-1}^h$, and, if $i \neq n$, then 
$x_i^h$ is preceded by $x_{i+1}^h$;

\smallskip
$\bullet$ if $i \neq n$, then $x_i^{-h}$ is followed by $x_{i+1}^{-h}$, and, if $i\neq 1$, then 
 $x_i^{-h}$ is preceded by $x_{i-1}^{-h}$.

\smallskip
Clearly, these properties hold if and only if $u$ is of the form $u=(x_n^h \dots x_2^h x_1^h)^r$ 
with $r \in \Z$.
\qed

\bigskip\noindent
{\bf Proof of Proposition A.1.} (1) Let $k \in \Z \setminus \{0\}$. Let $\phi: F_n \to F_n$ be 
the automorphism determined by $\phi(x_i)= x_i^{-1}$ for all $i=1, \dots, n$. One can easily 
verify that
$$
\phi^{-1} \circ \rho_k^{(1)} (\sigma_i) \circ \phi = \rho_{-k}^{(1)} (\sigma_i)
$$
for all $i=1, \dots, n-1$, thus $\rho_k$ and $\rho_{-k}$ are equivalent.

Let $k,l>0$. For a group $G$, we denote by $H_1(G)$ the abelianization of $G$, and, for a subgroup 
$H$ of $G$, we denote by $\langle\!\langle H \rangle\!\rangle$ the normal subgroup of $G$ 
generated by $H$. By Lemma A.3, we have
$$
F_n/ \langle\!\langle F_n^{ \rho_k^{(1)} (B_n)} \rangle\!\rangle \simeq \langle x_1, \dots, x_n \ 
|\ x_n^k \dots x_2^k x_1^k =1 \rangle,
$$
hence
$$
H_1( F_n/ \langle\!\langle F_n^{ \rho_k^{(1)} (B_n)} \rangle\!\rangle) \simeq (\Z/ k\Z) \times 
\Z^{n-1}.
$$
So, if $\rho_k^{(1)}$ and $\rho_l^{(1)}$ are equivalent, then $(\Z/ k\Z) \times \Z^{n-1} 
\simeq (\Z/ l\Z) \times \Z^{n-1}$, thus $k=l$.

\smallskip
(2) Write
$$
y_i= x_1^2 \dots x_{i-1}^2 x_i \quad {\rm for}\ i=1, \dots, n.
$$
One can easily verify that
$$
\rho^{(3)} (\sigma_k) (y_i) = \begin{cases}
y_i &\text{if }i\neq k,k+1\\
y_{k+1} &\text{if }i=k\\
y_{k+1} y_k^{-1} y_{k+1} &\text{if }i=k+1
\end{cases}
$$
Let $\phi: F_n \to F_n$ be the automorphism determined by $\phi(x_i) = y_{n-i+1}$ for $i=1, 
\dots, n$, and let $\mu: B_n \to B_n$ be the automorphism determined by $\mu(\sigma_i)= 
\sigma_{n-i}$ for $i=1, \dots, n-1$. From the expression of $\rho^{(3)} (\sigma_k) (y_i)$ given 
above, follows
$$
\phi^{-1} \circ \rho^{(3)}(\sigma_i)\circ \phi = \rho^{(2)} (\mu(\sigma_i))
$$
for all $i=1, \dots, n-1$, thus $\rho^{(2)}$ and $\rho^{(3)}$ are equivalent.

\smallskip
(3) Let $k>0$. For $u \in F_n$, we denote by $[u]$ the element of $H_1(F_n)\simeq\Z^n$ represented by 
$u$. We have
$$
\rho^{(2)} (\sigma_1^t) [x_1] = (t+1) [x_1] -t [x_2]
$$
for all $t \in \N$. On the other hand, $\rho_k^{(1)} (\beta)$ has finite order as an automorphism of 
$H_1(F_n)$, for all $\beta \in B_n$. This shows that $\rho^{(2)}$ and $\rho_k^{(1)}$ are not 
equivalent.
\qed


\bigskip\bigskip\noindent
\halign{#\hfill\cr
John Crisp\cr
Luis Paris\cr
Laboratoire de Topologie\cr
Universit\'e de Bourgogne\cr
UMR 5584 du CNRS, BP 47870\cr
21078 Dijon cedex\cr
FRANCE\cr
\noalign{\smallskip}
\texttt{ jcrisp@u-bourgogne.fr}\cr
\texttt{ lparis@u-bourgogne.fr}\cr}

\end{document}